\RequirePackage[english]{babel}
\documentclass{amsart}
\usepackage{amssymb,latexsym}
\usepackage{amsmath,amsfonts,amsthm}
\usepackage{stmaryrd}
\usepackage{enumerate}
\usepackage[T1]{fontenc}
\usepackage[latin1]{inputenc}
\usepackage[english]{babel}
\usepackage{mathtools}
\usepackage{url}

\usepackage[all]{xy}
\CompileMatrices
\SelectTips{cm}{}
\def\cxymatrix#1{\xy*[c]\xybox{\xymatrix#1}\endxy}

\usepackage{graphicx}
\usepackage{caption}
\usepackage{subfig}
\usepackage{graphicx}

\usepackage{ifpdf}
\ifpdf
\fi

\theoremstyle{plain}
\newtheorem{thm}{Theorem}[section]
\newtheorem{lemma}[thm]{Lemma}
\newtheorem{prop}[thm]{Proposition}
\newtheorem{cor}[thm]{Corollary}
\theoremstyle{definition}
\newtheorem{defn}[thm]{Definition}
\newtheorem{ex}[thm]{Example}
\newtheorem*{ack}{Acknowledgements}

\theoremstyle{remark}
\newtheorem{remark}[thm]{Remark}
\numberwithin{equation}{section}

\DeclareMathOperator{\id}{id}
\DeclareMathOperator{\Log}{Log}

\DeclareMathOperator{\Vol}{Vol}
\DeclareMathOperator{\CS}{CS}
\DeclareMathOperator{\cs}{cs}

\DeclareMathOperator{\SL}{SL}

\DeclareMathOperator{\PSL}{PSL}

\DeclareMathOperator{\Tor}{Tor}
\DeclareMathOperator{\Ker}{Ker}

\newcommand{\R}{\mathbb R} 
\newcommand{\C}{\mathbb C}
\renewcommand{\H}{\mathbb H}

\newcommand{\Cremove}{\C\backslash\{0,1\}}
\newcommand{\Mat}[4]{\left(\begin{smallmatrix}#1&#2\\#3&#4\end{smallmatrix}\right)}
 
\newcommand{\Z}{\mathbb Z}
\newcommand{\Q}{\mathbb Q}

\newcommand{\Ct}{\bar{C}}
\newcommand{\g}{\mathfrak g}
\newcommand{\B}{\mathcal B}
\newcommand{\Pre}{\mathcal P}

\begin{document}
\title[The complex volume of a representation]{The volume and Chern--Simons invariant\\of a representation}
\author{Christian~K. Zickert}
\address{Universiy of California\\Berkeley, CA 94720-3840, USA}
\email{zickert@math.berkeley.edu}
\begin{abstract} We give an efficient simplicial formula for the volume and Chern--Simons invariant of a boundary-parabolic $\PSL(2,\C)$--representation of a tame $3$--manifold. If the representation is the geometric representation of a hyperbolic $3$--manifold, our formula computes the volume and Chern--Simons invariant directly from an ideal triangulation with no use of additional combinatorial topology. In particular, the Chern--Simons invariant is computed just as easily as the volume.

\end{abstract} 
\maketitle

\section*{Introduction}
The volume and Chern--Simons invariant are two interesting and important invariants of a hyperbolic $3$--manifold. We will always assume that a hyperbolic $3$--manifold is complete, oriented and of finite volume, so that the hyperbolic structure is unique. Recall that the Chern--Simons invariant of a closed hyperbolic $3$--manifold $M$ is defined by the formula

\[\cs(M)=\frac{1}{8\pi^2}\int_{s(M)}\text{Tr}\bigl(A\wedge dA+\frac{2}{3}A\wedge A\wedge A\bigr)\in \R/\Z.\]
Here $A$ is the connection in the orthonormal frame bundle given by the hyperbolic metric, and $s(M)$ is an orthonormal frame field, i.e.~a section of the orthonormal frame bundle. The definition of the Chern--Simons invariant extends to hyperbolic manifolds with cusps using so-called ``special singular frame fields'' that are linear near the cusps. See Meyerhoff~\cite{Meyerhoff} for details. In the cusped case the Chern--Simons invariant is only defined modulo $1/2$. 

The Chern--Simons invariant is intimately related to the hyperbolic volume, and the two invariants are often regarded as the real and imaginary part of a so-called \emph{complex volume} given by
\begin{equation*}\Vol(M)+i\CS(M)\in \C/i\pi^2\Z,\end{equation*} where $\CS(M)=2\pi^2\cs(M)\in \R/\pi^2\Z$.

A very interesting feature of the complex volume is that it can be realized as a characteristic class of flat $\PSL(2,\C)$--bundles called the Cheeger--Chern--Simons class. 
This class satisfies that the characteristic cohomology class of the canonical flat $\PSL(2,\C)$--bundle over a closed hyperbolic $3$--manifold gives the complex volume when evaluated on the fundamental class.

The canonical flat $\PSL(2,\C)$--bundle over a hyperbolic $3$--manifold $M$ admits a unique bundle map to the universal flat $\PSL(2,\C)$--bundle over the classifying space $B(\PSL(2,\C))$, where $\PSL(2,\C)$ is regarded as a discrete group. This means that it is enough to study the characteristic class of the universal bundle. This class lies in $H^3(B(\PSL(2,\C)),\C/\pi^2\Z)$, and since $\C/\pi^2\Z$ is divisible, we can regard it as a homomorphism \[\hat c_2\colon H_3(\PSL(2,\C))\to \C/\pi^2\Z.\] 
As usual, we have identified the homology of a discrete group with the homology of its classifying space.
If $M$ is closed, the fundamental class of $M$ determines a fundamental class in $H_3(\PSL(2,\C))$, and the image of this class under the homomorphism $\hat c_2$ is $i(\Vol(M)+i\CS(M))$.

Note that if $M$ is closed, $\CS(M)$ is naturally defined in $\R/2\pi^2\Z$, whereas the image of $\hat c_2$ is in $\C/\pi^2\Z$. It is known that $\hat c_2$ does not admit a lift to $\C/2\pi^2\Z$, so the $2$--torsion of the Chern--Simons invariant of a closed hyperbolic manifold is not detected by the fundamental class.

In \cite{Neumann} Neumann has obtained an explicit formula for $\hat c_2$, and the computation of the complex volume of a closed hyperbolic manifold thus amounts to determining its fundamental class in $H_3(\PSL(2,\C))$. This, however, is quite difficult in general. Neumann gets around this by constructing a group $\widehat \B(\C)$, called the \emph{extended Bloch group}, which is isomorphic to $H_3(\PSL(2,\C))$, but more suitable for geometric purposes. He defines a map $R\colon \widehat \B(\C)\to \C/\pi^2\Z$ using an extended version of Rogers dilogarithm, and shows that under the identification of $\widehat \B(\C)$ with $H_3(\PSL(2,\C))$, the map $R$ corresponds to $\hat c_2$. He shows that every hyperbolic manifold defines an element in $\widehat \B(\C)$ whose image under $R$ is the complex volume (times $i$), and he thus obtains a formula for the complex volume that applies to cusped manifolds as well. 

Neumann's formula has been implemented in Snap, a freely available computer program~\cite{snapprogram} for numerical computation of invariants of hyperbolic $3$--mani\-folds. The formula works quite efficiently for manifolds with few simplices, but it involves some complicated combinatorial topology which slows down computations remarkably when the number of simplices increases. For example, if the number of simplices is around thirty or so, it will generally take Snap more than half an hour to compute the complex volume. 

In this paper we present a new approach which makes use of the relative homology group $H_3(\PSL(2,\C),P)$, where $P$ is the subgroup of upper triangular matrices with $1$ on the diagonal. We show in Section \ref{compoftrunc} that this group can be computed using a complex generated by ideal hyperbolic simplices endowed with a decoration consisting of a horosphere at each ideal vertex together with an identification of the horosphere with $\C$. Such a decoration naturally endows each ideal simplex with a flattening, and we use this to define a map 
\[\Psi\colon H_3(\PSL(2,\C),P)\to \widehat\B(\C).\] The formula is direct, and involves no combinatorial topology.

In section \ref{fundclasssec} we show that a tame $3$--manifold $M$ with a boundary-parabolic $\PSL(2,\C)$--representation $\rho$ defines a fundamental class in $H_3(\PSL(2,\C),P)$, which is defined once we have picked a decoration of $\rho$ consisting of a choice of conjugation of the $\rho$--image of each peripheral subgroup into $P$. 
Given a topological triangulation of $M$, we can construct an explicit representative of the fundamental class in the complex of decorated ideal simplices mentioned above. This is done using a developing map, whose purpose is to endow each simplex in the triangulation with the shape of an ideal simplex. The developing map naturally translates the decoration of $\rho$ into a decoration of the ideal simplices by horospheres. We stress that $M$ does not have to be hyperbolic, and the boundary components of $\bar M$ do not have to be tori.

The image of the fundamental class in $\widehat \B(\C)$ turns out to be independent of the choice of decoration, and we can define the complex volume of a boundary-parabolic representation $\rho$ by the formula \[i(\Vol(\rho)+i\CS(\rho))=R\circ\Psi ([\rho]), \] where $[\rho]$ is a fundamental class.
The formula agrees with that of Neumann in the special case where $\rho$ is the geometric representation of a hyperbolic $3$--manifold. Since every step of the process is natural and explicit, it allows us to compute the complex volume in an instant even for manifolds with a high number of simplices. As an example, we compute the complex volumes of all boundary-parabolic representations of the $5_2$ knot complement.

The set of boundary-parabolic representations of $M$ is often finite, and the set of complex volumes of these is an invariant of $M$. If $M$ is hyperbolic, this invariant can be viewed as a generalization of the Borel regulator of $M$, which consists of the set of volumes of the Galois conjugates of the geometric representation.


Most of the theory works in a more general setup. We show in Section \ref{fundclasssec} that any $G$--representation mapping boundary curves to conjugates of a fixed subgroup $H$, defines a fundamental class in $H_3(G,H)$, which is defined up to a choice of decoration. In the general setup, a decoration is a choice of element in the normalizer quotient $N_G(H)/H$ for each end. The generality of this approach suggests that the theory might have applications to the Chern--Simons theory of other Lie groups. 

Section \ref{SLreps} is a brief discussion of representations in $\SL(2,\C)$. We show that a cusped hyperbolic manifold with a spin structure determines a fundamental class in $H_3(\SL(2,\C))$ which is defined up to $2$--torsion. This $2$--torsion ambiguity is intrinsic, and has the interesting consequence that a large class of cusped hyperbolic manifolds, including hyperbolic knot complements, don't have ideal triangulations admitting strong, even valued flattenings. This may be true for all hyperbolic manifolds.
In Dupont--Zickert \cite{DupontZickert} and Goette--Zickert \cite{GZ} we obtained formulas for the Cheeger--Chern--Simons class $\hat c_2\colon H_3(\SL(2,\C))\to \C/4\pi^2\Z$, which is related to hyperbolic manifolds with spin structures. The formulas use even valued flattenings, and it is thus not clear how to apply these formulas to hyperbolic manifolds. In fact, the results of the present paper were derived in an attempt to improve this.

\begin{ack}
I wish to thank Walter Neumann for numerous enlightening discussions about this work and for his comments on preliminary versions of this paper. I also wish to thank Stavros Garoufalidis, Johan Dupont, Marc Culler and Charlie Frohmann for their interest in my work. Finally, I wish to thank the Danish grant ``Rejselegat for Matematikere'' for financial support.
\end{ack}

\section{The extended Bloch group}\label{extbsec} 
In this section we recall the definition of the extended Bloch group and some of its basic properties. The general reference for this is Neumann~\cite{Neumann}.
We start by recalling the definition of the classical Bloch group.

\begin{defn}\label{predefn} The \emph{pre-Bloch group} $\Pre(\C)$ is an abelian group generated by symbols $[z]$, $z\in \Cremove$ subject to the relation
\begin{equation}\label{fiveterm}
[x]-[y]+[\frac{y}{x}]-[\frac{1-x^{-1}}{1-y^{-1}}]+[\frac{1-x}{1-y}]=0.
\end{equation} 
This relation is called the \emph{five term relation}.
\end{defn}
\begin{defn}\label{blochdefn} The \emph{Bloch group} $\B(\C)$ is the kernel of the homomorphism
\[\nu\colon \Pre(\C)\to \wedge^2_\Z(\C^*)\]
defined by mapping a generator $[z]$ to $z\wedge (1-z)$. 
\end{defn}

Let $\H^3$ denote hyperbolic $3$--space and let $\bar \H^3$ denote its standard compactification. 
Unless otherwise specified, we shall always use the upper half space model for $\H^3$, which provides us with a natural identification of $\partial \bar \H^3$ with $\C\cup \{\infty\}$. In all of the following we identify the group of orientation preserving isometries of $\H^3$ with $\PSL(2,\C)$. The action of $\PSL(2,\C)$ on $\H^3$ extends uniquely to an action on $\bar \H^3$, with the action on $\partial \bar\H^3=\C\cup\{\infty\}$ being given by fractional linear transformations.

An \emph{ideal simplex} is a geodesic $3$--simplex whose vertices $z_0,z_1,z_2,z_3$ all lie in $\partial\bar \H^3=\C\cup \{\infty\}$. We consider the vertex ordering as part of the data defining an ideal simplex. It is well known that the orientation preserving congruence class of an ideal simplex is given by the \emph{cross-ratio}
\begin{equation}\label{cross-ratio}z=[z_0:z_1:z_2:z_3]=\frac{(z_0-z_3)(z_1-z_2)}{(z_0-z_2)(z_1-z_3)}\in \Cremove.\end{equation}
An ideal simplex is flat if and only if the cross-ratio is real, and if it is not flat, the orientation given by the vertex ordering agrees with the orientation inherited from $\H^3$ if and only if the cross-ratio has positive imaginary part.
Since an ideal simplex is determined up to congruence by its cross-ratio, we can regard the pre-Bloch group as being generated by (congruence classes of) ideal simplices. In this picture the five term relation is equivalent to the relation 
\begin{equation}\label{fivetermrel} \sum_{i=0}^4(-1)^i[z_0:\dots:\hat z_i:\dots:z_4]=0,\end{equation}
which implies that an element in $\Pre(\C)$ is invariant under $2$--$3$ moves and $1$--$4$ moves of the ideal simplices. See e.g.~Neumann~\cite{Neumann} for a description of these moves. 

It easily follows from \eqref{cross-ratio} that an even permutation of the $z_i$'s replaces $z$ by one of three so-called
\emph{cross-ratio parameters}.
 \[z,\qquad z'=\frac{1}{1-z},\qquad z''=1-\frac{1}{z}.\]

In the following we let $\Log$ denote a particular branch of logarithm that we fix once and for all. In concrete examples we will always use the principal branch having imaginary part in the interval $(-\pi,\pi]$.

\begin{defn}\label{logpardefn} Let $\Delta$ be an ideal simplex with cross-ratio $z$. A \emph{flattening} of $\Delta$ is a triple of complex numbers of the form
\begin{multline}(w_0,w_1,w_2)=\bigl(\Log z+p\pi i,-\Log (1-z)+q\pi i,\\
-\Log (z)+\Log(1-z)-p\pi i-q\pi i\bigr)\end{multline}
with $p,q\in\Z$. We call $w_0,w_1$ and $w_2$ \emph{log-parameters}. Up to multiples of $\pi i$, the log-parameters are logarithms of the cross-ratio parameters.
\end{defn} 

One can show that the set of flattened simplices has a natural structure as a Riemann surface with four components corresponding to the parities of $p$ and $q$. It is a $\Z\times \Z$ cover of $\Cremove$. We will not need this here.

\begin{remark}\label{pqlogremark}
Note that the log-parameters uniquely determine $z$. We can thus write a flattening as $[z;p,q]$. This notation, however, depends on the choice of logarithm branch.
\end{remark}

\begin{figure}[ht]
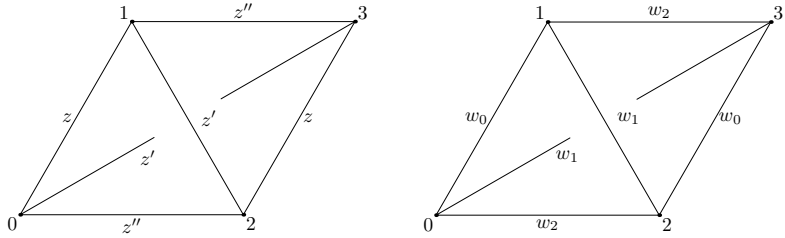

\centering
\begin{minipage}[c]{0.4\textwidth}
\centering
\includegraphics[width=4.8cm]{ChernSimonsFigure.8}
\end{minipage}
\quad
\begin{minipage}[c]{0.4\textwidth}
\centering
\includegraphics[width=4.8cm]{ChernSimonsFigure.9}
\end{minipage}
\captionsetup{width=9cm}
\caption{Associating cross-ratio parameters and log-parameters to edges of a flattened ideal simplex.}\label{asspar}
\end{figure}

In the following we will associate cross-ratio parameters and log-parameters to the edges of a flattened ideal simplex as indicated in Figure \ref{asspar}.

\begin{defn}\label{combflat}  
Let $z_0,\ldots,z_4$ be five distinct points in $\C\cup\{\infty\}$ and let $\Delta_i$ denote the simplices $[z_0,\dots,\hat z_i,\dots,z_4]$.
Suppose $(w_0^i,w_1^i,w_2^i)$ are flattenings of the simplices $\Delta_i$.
Every edge $[z_iz_j]$ belongs to exactly three of the $\Delta_i$'s 
and therefore has three associated log-parameters. 
The flattenings are said to satisfy the \emph{flattening condition} if for each edge the signed sum of the three associated log-parameters is zero. The sign is positive if and only if $i$ is even. 
\end{defn}
It follows directly from the definition that the flattening condition is equivalent to the following ten equations.
\begin{gather}\label{teneq}
\begin{aligned}
&[z_0z_1]:&w_0^2-w_0^3+w_0^4=0&\hspace{1cm} [z_0z_2]:&-w_0^1-w_2^3+w_2^4=0\\
&[z_1z_2]:&w_0^0-w_1^3+w_1^4=0&\hspace{1cm} [z_1z_3]:&w_2^0+w_1^2+w_2^4=0\\
&[z_2z_3]:&w_1^0-w_1^1+w_0^4=0&\hspace{1cm} [z_2z_4]:&w_2^0-w_2^1-w_0^3=0\\
&[z_3z_4]:&w_0^0-w_0^1+w_0^2=0&\hspace{1cm} [z_3z_0]:&-w_2^1+w_2^2+w_1^4=0\\
&[z_4z_0]:&-w_1^1+w_1^2-w_1^3=0&\hspace{1cm} [z_4z_1]:&w_1^0+w_2^2-w_2^3=0
\end{aligned}
\end{gather}

\begin{defn}\label{Phat} The \emph{extended pre-Bloch group} $\widehat \Pre(\C)$ is the free abelian group on flattened ideal simplices subject to the relations
\begin{enumerate}[(i)]
\item $\sum_{i=0}^4(-1)^i(w_0^i,w_1^i,w_2^i)=0$ if the flattenings satisfy the flattening condition.
\item $[z;p,q]+[z;p',q']=[z;p,q']+[z;p',q].$
\end{enumerate}
\end{defn}
The first relation obviously lifts the relation \eqref{fivetermrel}. It is therefore called the \emph{lifted five term relation}. The second relation is called the \emph{transfer relation} and it plays a more subtle role. We refer to Goette--Zickert~\cite{GZ} and Neumann~\cite{Neumann} for a discussion.
\begin{defn}\label{Bhat} The \emph{extended Bloch group} $\widehat \B(\C)$ is the kernel of the homomorphism
\[\widehat \nu\colon \widehat \Pre(\C)\to \wedge^2_\Z(\C) \]
defined on generators by $(w_0,w_1,w_2)\mapsto w_0\wedge w_1$.
\end{defn}
The relationship between the extended Bloch group and the classical Bloch group is summarized in the diagram below, which is taken from Neumann~\cite{Neumann}.
\begin{thm} There is a commutative diagram with exact rows and columns.
\[\xymatrix{&0\ar[d]&0\ar[d]&0\ar[d]&&\\0\ar[r]&\mu^*\ar[r]\ar[d]^-\chi&{\C^*}\ar[r]\ar[d]^\chi&{\C^*/\mu^*}\ar[r]\ar[d]^\beta&0\ar[d]&\\
0\ar[r]&{\widehat \B(\C)}\ar[r]\ar[d]^-\pi&{\widehat \Pre(\C)}\ar[r]^-{\widehat \nu}\ar[d]^-\pi&{\wedge^2_\Z(\C)}\ar[r]\ar[d]^\epsilon&{K_2(\C)}\ar@2{-}[d]\ar[r]&0\\0\ar[r]&{\B(\C)}\ar[r]\ar[d]&{\Pre(\C)}\ar[r]^-\nu\ar[d]&{\wedge^2_\Z(\C^*)\ar[d]}\ar[r]&{K_2(\C)}\ar[r]\ar[d]&0\\&0&0&0&0}\]
Here $\mu^*$ is the group of roots of unity, and the maps not already defined above, are defined as follows:
\begin{align*}
\chi(z)&=[z;0,1]-[z;0,0];\\
\beta([z])&=\textnormal{log}(z)\wedge \pi i;\\
\epsilon(x\wedge y)&=-\exp(x)\wedge \exp(y);\\
\pi([z;p,q])&=[z].
\end{align*}
\end{thm}
In \cite{Neumann} Neumann shows that the map 
\begin{gather}\widehat L\colon \widehat\Pre(\C)\to \C/\pi^2\Z\nonumber\\
[z;p,q]\mapsto L(z)+\frac{\pi i}{2} (q\Log(z)+p\Log(1-z))-\pi^2/6\label{lhatdefn}
\end{gather}
is well defined.
Here $L(z)=-\int_0^z\frac{\Log(1-t)}{t}dt+\frac{1}{2}\Log(z)\Log(1-z)$ is Rogers dilogarithm. The map $\widehat L$ is denoted by $R$ in Neumann \cite{Neumann}. 

\begin{remark} As mentioned in Remark \ref{pqlogremark}, the representation of a flattening as $[z;p,q]$ depends on the choice of logarithm branch, but one can check that the expression \eqref{lhatdefn} is independent of this choice.\end{remark}

\section{Relative homology of groups}\label{relative}
Recall that the homology of a discrete group $G$ is equal to the singular homology of its classifying space $BG$, and can be calculated as $H_*(F_*\otimes_G\Z)$, where $F_*$ is any free $G$--resolution of $\Z$.

Let $H$ be a subgroup of $G$ and let $Cof(i)$ denote the cofiber (mapping cone) of the map $BH\to BG$ induced by inclusion. We define the \emph{relative homology}, denoted $H_*(G,H)$, to be the reduced singular homology groups $\tilde H_*(Cof(i);\Z)$. Regarding $BH$ as a subspace of $BG$, this is isomorphic to $H_*(BG,BH;\Z)$. 

For any set $X$ we can construct a complex $C_*(X)$ of abelian groups by letting $C_n(X)$ be the free abelian group generated by $(n+1)$--tuples of elements of $X$. The boundary map is given by 
\[\partial(x_0,\dots,x_n)=\sum_{i=0}^n(-1)^i(x_0,\dots,\hat x_i,\dots,x_n).\]
The complex $C_*(X)$ is acyclic in dimensions greater than $0$ and $H_0(C_*(X))=\Z$.
If $X$ is a group $G$, left multiplication endows $C_n(G)$ with the structure of a free $G$--module and $C_*(G)$ becomes a free $G$--resolution of $\Z$. Hence, the complex  
\begin{equation}\label{defofB} 
B_*(G)=C_*(G)\otimes_{\Z[G]}\Z 
\end{equation}
calculates the homology of $G$.
Theorem~\ref{relhom} below gives a similar description of relative homology in terms of free resolutions. This is probably well known, but since we don't know of any reference we include a proof.

\begin{thm}\label{relhom}Let $H$ be a subgroup of $G$ and let $K$ be the kernel of the augmentation map $C_0(G/H)\to \Z$. For any free $G$--resolution $\{F_i\}_{i=1}^\infty$ of $K$ we have a canonical isomorphism \[H_*(F_*\otimes_{\Z[G]}\Z)\cong H_*(G,H).\]
That is, $H_*(G,H)=\Tor^{\Z[G]}_*(K,\Z)$.
\end{thm}
\begin{proof}
It is enough to prove the existence of a free $G$--resolution of $K$ for which the isomorphism holds.
Let $B_*(H)$ and $B_*(G)$ be as in \eqref{defofB} and let $i_*$ denote the map induced by inclusion. By the standard cone construction (see e.g.~Chapter $4.2$ in Spanier~\cite{Spanier}), the reduced homology of $Cof(i)$ is the homology of the complex $D_i=B_{i-1}(H)\oplus B_i(G)$, with boundary map given by the matrix $\Mat{\partial}{0}{i_*}{-\partial}$. Define a complex $F_i$ of free $G$--modules by 
\begin{align*} F_i&=D_i\otimes_\Z\Z[G],\quad \text{for }i\geq 2\\
 F_1&=\Ker (D_1\to D_0)\otimes_\Z\Z[G]=B_1(G)\otimes_\Z\Z[G]=C_1(G).
 \end{align*}
Note that $H_i(F_*\otimes_{\Z[G]}\Z)=H_i(D_*)$ for all $i\geq 1$. The theorem will now follow if we can prove that the map $F_2\to F_1$ has cokernel isomorphic to $K$. 

Define a map $\rho\colon F_1\to C_0(G/H)$ as the composition
\[\xymatrix{F_1=C_1(G)\ar[r]^-\partial&C_0(G)\ar[r]^-\pi&C_0(G/H),}\]
where $\pi$ is induced by projection onto cosets. It is now simple to check that $\rho$ maps surjectively onto $K$ with kernel equal to the image of $F_2$. This proves the theorem.
\end{proof}

The complex $C_*(G)$ can be regarded as being generated by simplices with a $G$--labeling of vertices. Also, $B_*(G)$ can be regarded as being generated by simplices with a $G$--labeling of edges. See e.g.~Chapter IV in Mac Lane~\cite{MacLane} for more explanation. Similar to this, relative homology can be computed using complexes of \emph{truncated} simplices with labelings. This will be explored in the next section.

\section{The complex of truncated simplices}\label{compoftrunc}
Let $G=\PSL(2,\C)$ and let $P$ be the image in $G$ of the group of upper triangular matrices with $1$ on the diagonal. Note that $P$ is isomorphic to $\C$. We now construct an explicit complex computing the relative homology groups $H_*(G,P)$.

Let $\Delta$ be an $n$--simplex with a vertex ordering given by associating an integer $i\in\{0,\dots,n\}$ to each vertex. Let $\bar \Delta$ denote the corresponding truncated simplex obtained by chopping off disjoint regular neighborhoods of the vertices. Each vertex of $\bar \Delta$ is naturally associated with an ordered pair $ij$ of distinct integers. Namely, the $ij$--th vertex of $\bar \Delta$ is the vertex near the $i$--th vertex of $\Delta$ and on the edge going to the $j$--th vertex of $\Delta$.
\begin{defn}\label{trunccomp} Let $\Ct_n(G,P)$, $n\geq 1$, be the free abelian group generated by $G$--labelings $\{g^{ij}\}$ of vertices of truncated $n$--simplices satisfying
\begin{enumerate}[(i)]\item\label{one} For fixed $i$ the vertices $ij$ are labeled by distinct elements in $G$ mapping to the same left $P$--coset.
\item\label{two} The elements $g_{ij}=(g^{ij})^{-1}g^{ji}$ are counter diagonal, i.e.~of the form $\Mat{0}{-a^{-1}}{a}{0}$.
\end{enumerate}\end{defn}

Left multiplication endows $\Ct_n(G,P)$ with a $G$--module structure and the usual boundary map induces a boundary map on $\Ct_*(G,P)$ making it into a chain complex.

\begin{remark}
We will prove later that $\bar C_*(G,P)\otimes_{\Z[G]}\Z$ computes the relative homology groups $H_*(G,P)$.
For this to hold, property \eqref{two} of Definition \ref{trunccomp} is not required. Nor is distinctness in Property \eqref{one}. In fact, we can define $\bar C_*(G,H)$ for an arbitrary group $G$ and an arbitrary subgroup $H$ exactly as in Definition \ref{trunccomp} but without Property \eqref{two} and without distinctness in Property \eqref{one}. The equality $H_*(G,H)=H_*(\bar C_*(G,H)\otimes_{\Z[G]}\Z)$ will still hold. The reason for adding the extra properties is that we will be able to interpret a generator as an ideal simplex which is naturally flattened. This will be explained in the next section. 
\end{remark}

\begin{figure}[ht]
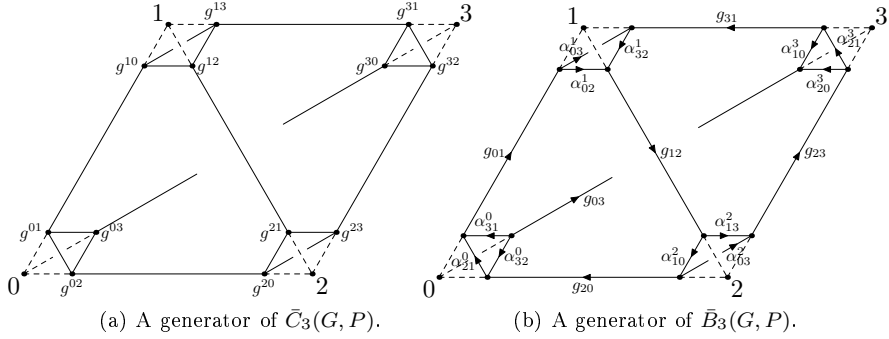

\centering
\subfloat[\label{sub}A generator of $\bar C_3(G,P)$.]
{\includegraphics[width=6.2cm]{ChernSimonsFigure.1}}\hspace{-0.8cm}
\subfloat[\label{sub2}A generator of $\bar B_3(G,P)$.]
{\includegraphics[width=6.2cm]{ChernSimonsFigure.13}}
\captionsetup{width=9cm}
\caption{Generators in the truncated complexes. In (b), the product of labelings around each two-face is $1$. For the front face, this implies that $g_{01}\alpha^1_{02}g_{12}\alpha^2_{10}g_{20}\alpha^0_{21}=1$. Reversing an arrow corresponds to replacing the appropriate label by its inverse.}\label{main}
\end{figure}

Note that $\Ct_n(G,P)$ is a free $G$--module, and that we can represent a generator by a truncated simplex together with a labeling of each oriented edge, such that an edge going from vertex $ij$ to vertex $kl$ is labeled by $(g^{ij})^{-1}g^{kl}$. We denote the labeling of an edge going from vertex $i$ to $j$ in the untruncated simplex by $g_{ij}$, and the labeling of the edges near the $k$--th vertex by $\alpha^k_{ij}$, see Figure \ref{main}. We call these edges the \emph{long} edges and the \emph{short} edges, respectively. By properties \eqref{one} and \eqref{two} of Definition \ref{trunccomp}, the $\alpha^k_{ij}$'s are non-trivial elements in $P$ and the $g_{ij}$'s are counter diagonal. Furthermore, the edge labelings are forced to satisfy that the product of labelings along any two-face (including the triangles) is $1$. We denote the complex generated by such $G$--labelings $\bar B_*(G,P)$. By definition we have $\bar B_*(G,P)=\bar C_*(G,P)\otimes_{\Z[G]} \Z$.

In the following we will often regard the labelings of short edges as complex numbers using the canonical identification of a matrix $\Mat{1}{x}{0}{1}\in P$ with the complex number $x\in \C$. The matrix corresponding to $x$ will be denoted $(x)$.

Note that the labelings of the short edges (regarded as complex numbers) satisfy

\begin{equation}
\begin{gathered}\label{gather}
\alpha^i_{jk}=-\alpha^i_{kj},\\
\alpha^i_{jk}+\alpha^i_{kl}+\alpha^i_{lj}=0.
\end{gathered}
\end{equation}

\begin{lemma}\label{unique} Let $\alpha\in \bar B_n(G,P)$ be a generator. For any $i,j,k,l$ the labelings of the short edges (regarded as complex numbers) satisfy
\begin{equation}\label{alphas}
\alpha^i_{kj}\alpha^j_{ik}=\alpha^i_{lj}\alpha^j_{il}.
\end{equation}
Moreover, if $\alpha^i_{jk}\neq 0$ are labelings of the short edges of a truncated $n$--simplex satisfying \eqref{alphas} and \eqref{gather}, there is a unique way of labeling the long edges to obtain a generator of $\bar B_n(G,P)$.  
\end{lemma}

\begin{proof} Consider the two-face of $\alpha$ defined by the vertices $i$, $j$ and $k$. For notational simplicity, assume that the labelings of the edges are given by
\begin{gather*}
g_{ij}=\Mat{0}{-a^{-1}}{a}{0},\quad g_{jk}=\Mat{0}{-b^{-1}}{b}{0},\quad g_{ki}=\Mat{0}{-c^{-1}}{c}{0},\\
\big(\alpha^j_{ik}\big)=\Mat{1}{p}{0}{1},\quad \big(\alpha^k_{ji}\big)=\Mat{1}{q}{0}{1},\quad \big(\alpha^i_{kj}\big)=\Mat{1}{r}{0}{1}.\end{gather*}
A simple calculation shows that \[g_{ij}\big(\alpha^j_{ik})g_{jk}\big(\alpha^k_{ji})g_{ki}\big(\alpha^i_{kj})=\Mat{\frac{-bc}{a}q}{\frac{-b}{ac}(c^2qr-1)}{\frac{ac}{b}(b^2pq-1)}{\frac{a}{bc}\big((b^2pq-1)c^2r-b^2p\big)},\] which is easily seen to be $1\in \PSL(2,\C)$ if and only if 
\begin{equation}\label{squareroot}a^{-2}\!=rp=\alpha^i_{kj}\alpha^j_{ik},\quad b^{-2}\!=pq=\alpha^j_{ik}\alpha^k_{ji},\quad c^{-2}\!=qr=\alpha^k_{ji}\alpha^i_{kj}.\end{equation}
Since $a$ only depends on $i$ and $j$, the first statement follows. Given labelings of the short edges, we can use \eqref{squareroot} to define the labelings of the long edges. This is consistent since the $\alpha_{ij}^k$'s satisfy \eqref{alphas}.
\end{proof}
\begin{remark}\label{SL} It follows from equation~\eqref{squareroot} that the square root defined by our choice of logarithm branch gives us particular representatives of the $g_{ij}$'s in $\SL(2,\C)$ satisfying $g_{ij}=g_{ji}$. In the following we shall thus always regard the $g_{ij}$'s as elements in $\SL(2,\C)$. Note, however, that the product along hexagonal faces may now be $-1$ instead of $1$ in $\SL(2,\C)$.
\end{remark}
We conclude the section with a proof that $\bar B_*(G,P)$ computes the relative homology groups $H_*(G,P)$. 

\begin{lemma}\label{uniquerep} Let $gP$ and $hP$ be $P$--cosets, satisfying that $gB\neq hB$, where $B=\{\Mat{\lambda}{z}{0}{\lambda^{-1}}\in \PSL(2,\C)\}$. There exist unique coset representatives $gx$ and $hy$ satisfying that $(gx)^{-1}hy$ is counter diagonal.
\end{lemma} 
\begin{proof}
Let $g^{-1}h=\Mat{a}{b}{c}{d}$ and let $x=\Mat{1}{p}{0}{1}$ and $y=\Mat{1}{q}{0}{1}$.
We have 
\[x^{-1}g^{-1}hy=\Mat{a-cp\,}{\,aq+b-p(cq+d)}{c}{cq+d}.\]
Since $gB\neq hB$, it follows that $c$ is non-zero. This implies that there exist unique complex numbers $p$ and $q$ such that the above matrix is counter diagonal. The elements $gx$ and $hy$ are easily seen to be independent of the representatives $g$ and~$h$. 
\end{proof}

\begin{cor}\label{cosetdec} 
An $(n\!+\!1)$--tuple of left $P$--cosets that are distinct as cosets of $B$ uniquely determines a generator of $\bar C_n(G,P)$.
\end{cor}

\begin{prop}\label{freeres} Let $K$ be the kernel of the augmentation map $C_0(G/P)\to \Z$. The complex $\Ct_*(G,P)$ gives a free $G$--resolution of $K$.
\end{prop}
\begin{proof}  We have already seen that $\Ct_*(G,P)$ is a complex of free modules. To prove exactness we use a standard cone argument:
Given a generator $\alpha\in \Ct_n(G,P)$, we wish to define a ``cone'' $S(\alpha)\in \bar C_{n+1}(G,P)$. Let $h_iP$ be the coset determined by vertex $i$ of $\alpha$ and let $gP$ be a coset satisfying that $gB\neq h_iB$ for all $i$. The cone will depend on this choice of coset.
Let $\bar \Delta^{n+1}$ be a truncated $(n+1)$--simplex. For $i,j\neq 0$ define the labeling of vertex $ij$ to be the labeling of vertex $(i-1)(j-1)$ of $\alpha$. The remaining vertices can be labeled using Lemma~\ref{uniquerep} above. Namely, we label vertex $0i$ by $gx_i$ and vertex $i0$ by $h_{i-1}y_i$, where $x_i$ and $y_i$ are defined as in Lemma~\ref{uniquerep}. This finishes the definition of $S(\alpha)$. We can similarly define a cone on any chain $\beta\in \bar C_n(G,P)$ as long as $gB$ is distinct from the $B$--cosets determined by the summands. Since $G/B$ is infinite we can always find such a coset. It follows directly from the construction that 
 \[\partial (S(\beta))-S(\partial(\beta))=\beta.\]
This shows that every cycle is a boundary.
The only thing left to prove is that the map $\partial_2\colon \Ct_2(G,P)\to \Ct_1(G,P)$ has cokernel isomorphic to $K$. 
Let $\pi\colon \bar C_1(G,P)\to C_0(G/P)$ be the map induced by $\partial_1$. It is trivial to see that $\pi$ has image in $K$ and maps the image of $\partial_2$ to $0$. 
We need to prove that each chain in the kernel of $\pi$ lies in the image of $\partial_2$. For a generator $\alpha\in \Ct_1(G,P)$ we write $\partial_1(\alpha)=\alpha^1-\alpha^0$, with $\alpha^i$ generators of $C_0(G)$. Since the complex $C_*(G/P)$ is acyclic, we can write any chain in the kernel of $\pi$ as a sum of chains of the form $\alpha_0-\alpha_1+\alpha_2$ satisfying that $\alpha_i^1$ and $\alpha_{i+1}^0$ (indices modulo $3$) are in the same $P$--coset.
Define $g^{ij}$ by the formula (indices modulo $3$)
\[\partial \alpha_i=\alpha^1_i-\alpha^0_i=g^{i+2,i+1}-g^{i+1,i+2},\] 
and let these be the labelings of a generator $\tau\in\Ct_2(G,P)$. It is now a simple matter to check that $\partial_2 \tau=\alpha_0-\alpha_1+\alpha_2$. This concludes the proof.
\end{proof} 
\begin{cor}\label{Hhom} We have an isomorphism 
\[H_*(\bar B_*(G,P))=H_*(G,P).\]
\end{cor}
\begin{proof}
This follows immediately from Theorem \ref{relhom}.
\end{proof}

\subsection{Decorations and flattenings}\label{crsection}
In this section we discuss some of the underlying geometry behind the complex of truncated simplices.
We shall see that every generator of $\bar B_3(G,P)$ can be regarded as an ideal simplex together with a decoration which endows the ideal simplex with a natural flattening. This will be used to define a map $\Psi\colon H_3(G,P)\to \widehat \B(\C)$. 

Recall that $G$ acts on $\H^3$, which we identify with upper half space. The subgroup $P$ fixes $\infty\in \bar \H^3$ and acts by translations on any horosphere at $\infty$. We endow a horosphere at $\infty$ with the counterclockwise orientation as viewed from $\infty$. Using the action of $G$, this induces an orientation on all horospheres.
\begin{defn}\label{euclidhorodef} A horosphere together with a choice of orientation preserving isometry to $\C$ is called a \emph{euclidean horosphere}. We consider two euclidean horospheres based at the same point equal if the isometries differ by a translation. We let $G$ act on the set of euclidean horospheres in the obvious way.
\end{defn}

\begin{remark}\label{Hinf}
The action of $G$ on the set of euclidean horospheres is transitive with stabilizer $P$. Hence, the set of euclidean horospheres can be identified with the set of left $P$-cosets. This identification is fixed once we have picked a reference euclidean horosphere. A natural choice for such is the plane (horosphere at $\infty$) at height $1$ over the $x$--$y$ plane (identified with $\C$) with the euclidean structure induced by projection. For future reference we will denote this by $H(\infty)$.
\end{remark} 

\begin{defn}\label{truncsimp}
A choice of euclidean horosphere at each vertex of an ideal simplex is called a \emph{decoration} of the simplex. Having fixed a decoration, we say that the ideal simplex is \emph{decorated}. Two decorated ideal simplices are called \emph{congruent} if they differ by an element in $G$.
\end{defn}

\begin{figure}[ht]
\centering
\includegraphics[width=4cm]{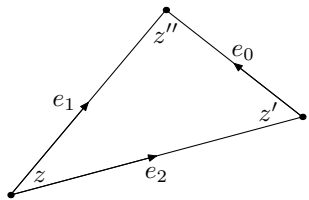}
\captionsetup{width=9cm}
\caption{An intersection triangle of an ideal simplex. The cross-ratio parameter at a vertex indicates the relationship between the two outgoing edges. Regarding the oriented edges as complex numbers, we have $e_1=ze_2$, $e_0=-z'^{-1}e_2$ and $e_0=z''e_1$. 
The ordering of the edges is the one induced from the vertex ordering of the ideal simplex.}\label{Ctriangle}
\end{figure}

A horosphere based at one of the ideal vertices of an ideal simplex intersects the simplex in an oriented euclidean triangle, which we will refer to as an \emph{intersection} triangle. 
A decoration enables us to view the intersection triangles as explicit triangles in $\C$. The association of cross-ratio parameters to the edges of an ideal simplex (see Figure~\ref{asspar}) associates cross-ratio parameters to the vertices of the intersection triangles as shown in Figure~\ref{Ctriangle}.

We wish to show that there is a one-to-one correspondence between generators of $\bar B_3(G,P)$ and congruence classes of decorated ideal simplices. 
Let $\{g^{ij}\}$ be vertex labelings of a generator in $\Ct_3(G,P)$ corresponding to an inhomogeneous generator $\alpha\in\bar B_3(G,P)$. Define ideal vertices $v_i\in\C\cup\{\infty\}$ by $v_i=g^{ij}\infty$. This is independent of $j$ by property \eqref{one} of Definition \ref{trunccomp}. The $v_i$'s determine an ideal simplex $\Delta$, which up to congruence only depends on $\alpha$. 
Since the product of edge labelings of $\alpha$ along each cut-off triangle is $1$, the identification of $P$ with $\C$ gives each cut-off triangle the geometric shape of a euclidean triangle determined up to translation.
We wish to prove that the cut-off triangles correspond to intersection triangles of a decoration of $\Delta$. A necessary condition for this to hold is that the cut-off triangles satisfy the same geometric properties as the geometric properties of the intersection triangles described above. 

\begin{lemma}\label{simclass} Let $\alpha\in \bar B_3(G,P)$ be a generator and let $z$ denote the cross-ratio of the associated ideal simplex $\Delta$. The labelings satisfy the following relations. In particular, the cut-off triangles are similar.
\begin{gather}\label{alphacr}
\begin{aligned}
z=\alpha^0_{12}/\alpha^0_{13}=\alpha^1_{03}/\alpha^1_{02}=\alpha^2_{30}/\alpha^2_{31}=\alpha^3_{21}/\alpha^3_{20}\\
z'=\alpha^0_{13}/\alpha^0_{23}=\alpha^3_{02}/\alpha^3_{01}=\alpha^1_{20}/\alpha^1_{23}=\alpha^2_{31}/\alpha^2_{01}\\
z''=\alpha^0_{32}/\alpha^0_{12}=\alpha^2_{01}/\alpha^2_{03}=\alpha^1_{23}/\alpha^1_{03}=\alpha^3_{10}/\alpha^3_{12}
\end{aligned}
\end{gather}
\end{lemma}
\begin{proof}
We only need to prove the first equation. The other two are trivial consequences.
Recall that the vertices of $\Delta$ are given by a choice of homogeneous representative of $\alpha$ in $\bar C_3(G,P)$. 
Using the unique representative with $g^{01}\!=\!1$ in $\PSL(2,\C)$, we obtain that the ideal vertices are $v_0=\infty$, $v_1=g_{01}\infty$, $v_2=(\alpha_{12}^0)g_{02}\infty$ and $v_3=(\alpha_{13}^0)g_{03}\infty$, and computing the cross-ratio yields \[[\infty,g_{01}\infty,(\alpha_{12}^0)g_{02}\infty,(\alpha_{13}^0)g_{03}\infty]=[\infty,0,\alpha_{12}^0,\alpha_{13}^0]=\alpha^0_{12}/\alpha^0_{13}=z.\]
Doing the same for the representative with $g^{23}=1$ gives vertices $v_0=\alpha^2_{30}$, $v_1=\alpha^2_{31}$, $v_2=\infty$, $v_3=0$ and cross-ratio $\alpha^2_{30}/\alpha^2_{31}$. This proves the first and third equality. The second and fourth equality follow directly from \eqref{alphas}.
\end{proof}

\begin{thm}\label{trunconetoone} Generators of $\bar B_3(G,P)$ are in one-one cor\-respon\-dence with congruence classes of decorated ideal simplices.
\end{thm} 
\begin{proof}
We have already seen that the decoration of a decorated ideal simplex endows each intersection triangle with the shape of a triangle in $\C$ which is determined up to translation. Using the identification of complex numbers with elements in $P$, this determines labelings of the small edges in a truncated simplex. Note that these labelings only depend on the congruence class of the decorated simplex. Using the geometry of the intersection triangles described in Figure~\ref{Ctriangle}, we see that the labelings satisfy \eqref{alphacr}. Since \eqref{alphacr} obviously implies \eqref{alphas} it follows from Lemma~\ref{unique} that the labelings define a unique generator of $\bar B_3(G,P)$. Now let $\alpha\in \bar B_3(G,P)$ be a generator and let $\Delta$ be the associated ideal simplex. Since the cut-off triangles satisfy \eqref{alphacr}, there is a unique way of picking euclidean horospheres at the ideal vertices of $\Delta$ such that the intersection triangles coincide with the cut-off triangles.
\end{proof}

\begin{remark}\label{newviewpoint}
Theorem \ref{trunconetoone} also follows from Remark \ref{Hinf} and Corollary \ref{cosetdec}. Note that if we identify the euclidean horosphere at a vertex of a decorated ideal simplex with a coset $gP$, then $g^{-1}$ takes the intersection triangle to a triangle in $H(\infty)$ whose projection onto $\C$ equals (up to translation) the explicit triangle given by the decoration.
\end{remark}

As we shall see below, a decorated ideal simplex is naturally equipped with a flattening. 
For a matrix $g=\Mat{a}{b}{c}{d}$ we let $c(g)$ denote the entry $c$. 
\begin{lemma}\label{flattening} Let $\alpha$ be a generator of $\bar B_3(G,P)$, and let $z$ denote the cross-ratio of $\alpha$ regarded as an ideal simplex. Note that for each of the long edges, $c(g_{ij})$ is well defined by Remark \ref{SL} and is non-zero. We have
\begin{equation}\label{crsigns}\frac{c(g_{03})c(g_{12})}{c(g_{02})c(g_{13})}=\pm z,\quad \frac{c(g_{13})c(g_{02})}{c(g_{01})c(g_{23})}=\pm \frac{1}{1-z},\quad \frac{c(g_{01})c(g_{23})}{c(g_{03})c(g_{12})}=\pm (1-\frac{1}{z}).\end{equation}
\end{lemma}
\begin{proof}  
From \eqref{squareroot}, we have that $c(g_{ij})^2=(\alpha^i_{kj}\alpha^j_{ik})^{-1}=(\alpha^i_{lj}\alpha^j_{il})^{-1}$ and using \eqref{alphacr} one easily checks that 
\[\frac{c(g_{03})^2c(g_{12})^2}{c(g_{02})^2c(g_{13})^2}=z^2,\]
from which the first equality follows. The other equalities are proved similarly. 
\end{proof}

\begin{remark}\label{pqeven}
One can prove that the signs in \eqref{crsigns} are $(+,+,-)$ if and only if the product of labelings around hexagonal faces is constant. By Remark \ref{SL} this constant is either $1$, in which case $\alpha$ is in $\bar B_3(\SL(2,\C),P)$, or $-1$, in which case $-\alpha$ is in $\bar B_3(\SL(2,\C),P)$. Here $-\alpha$ denotes the labeled truncated simplex obtained from $\alpha$ by changing the signs of all the labelings of long edges. We will need this in Section \ref{SLreps} where we discuss $\SL(2,\C)$--representations and their relations to even flattenings.
\end{remark}

It now follows from Lemma~\ref{flattening} that we can define a flattening of $\Delta$ by defining log-parameters
\begin{gather}\label{logpar}
\begin{aligned}
w_0(\alpha)=\Log c(g_{03})+\Log c(g_{12})-\Log c(g_{02})-\Log c(g_{13}),\\
w_1(\alpha)=\Log c(g_{02})+\Log c(g_{13})-\Log c(g_{01})-\Log c(g_{23}),\\
w_2(\alpha)=\Log c(g_{01})+\Log c(g_{23})-\Log c(g_{03})-\Log c(g_{12}).
\end{aligned}
\end{gather}
Note that since $c(g_{ij})=c(g_{ji})$, the log-parameters above are a sum of logarithms of complex numbers associated to the unoriented edges of $\Delta$. We will refer to them as $\Log(c)$--parameters.

Consider the map 
\begin{equation} \Psi\colon\bar B_3(G,P)\to \widehat\Pre(\C),\quad \alpha \mapsto (w_0(\alpha),w_1(\alpha),w_2(\alpha)).\end{equation}
\begin{thm}
The map $\Psi$ defined above sends boundaries to $0$ and cycles to $\widehat \B(\C)$, and we therefore obtain an induced map  
\begin{equation}\label{PsiHom}
\Psi\colon H_3(G,P)\to \widehat \B(\C).
\end{equation}
\end{thm}
\begin{proof}
This follows as in Section $3.1$ in Dupont--Zickert~\cite{DupontZickert}.
Let $\alpha\in \bar B_4(G,P)$ be a generator, and let $\alpha_i$ denote the five generators of $\bar B_3(G,P)$ obtained from $\alpha$ by deleting the $i$--th vertex. We have \[\partial\alpha=\sum (-1)^i\alpha_i.\]
Let $(w_0^i,w_1^i,w_2^i)$ be the flattening of the simplex corresponding to $\alpha_i$ as defined by \eqref{logpar}. 
We have that $\Psi(\partial\alpha)=0$ if and only if the flattenings $(w_0^i,w_1^i,w_2^i)$ satisfy the flattening condition, which is equivalent to satisfying the ten equations \eqref{teneq}.
We verify the first of these and leave the others to the reader.
We have 
\begin{gather*}
\begin{aligned}
w_0^2=&\Log c(g_{04})+\Log c(g_{13})-\Log c(g_{03})-\Log c(g_{14})\\
w_0^3=&\Log c(g_{04})+\Log c(g_{12})-\Log c(g_{02})-\Log c(g_{14})\\
w_0^4=&\Log c(g_{03})+\Log c(g_{12})-\Log c(g_{02})-\Log c(g_{13}),
\end{aligned}
\end{gather*} 
from which it follows that the equation $w_0^2-w_0^3+w_0^4=0$ is satisfied.
Having verified all the ten equations of \eqref{teneq}, we have proved that $\Psi$ sends boundaries to zero.
To prove that $\Psi$ sends cycles to $\widehat \B(\C)$ define a map $\mu\colon \bar B_2(G,P)\to \C\wedge_\Z \C$ by
\begin{multline} \alpha\mapsto \Log c(g_{01})\wedge\Log c(g_{02})\\-\Log c(g_{01})\wedge\Log c(g_{12})+\Log c(g_{02})\wedge\Log c(g_{12}).
\end{multline}
Letting $\Z[\widehat\C]$ be the free abelian group on the set of flattenings, a straightforward calculation shows that the diagram below is commutative.
\begin{equation}\cxymatrix{{\bar B_3(G,P)\ar[d]^\partial\ar[r]^-\Psi&{\Z[\widehat\C]}\ar[d]^{\widehat\nu}\\\bar B_2(G,P)\ar[r]^-\mu&{\C\wedge_\Z\C}}}\end{equation}
This means that cycles are mapped to $\widehat \B(\C)$ as desired.
\end{proof}
\begin{remark}\label{depsqrt} The reader may argue that $\Psi$ depends on the choice of logarithm used to define the flattenings. This is not the case, see Remark \ref{choice}.\end{remark}


The composition $\widehat L\circ \Psi\colon H_3(G,P)\to \C/\pi^2\Z$, where $\widehat L$ is given by \eqref{lhatdefn}, can be viewed as a relative Cheeger--Chern--Simons class. The fact that it agrees with $\hat c_2$ on $H_3(G)$ will be proved in Section~\ref{compvolsec}.

\section{Boundary-parabolic representations}\label{parsec}
In this section we define the notion of a boundary-parabolic representation of a tame manifold and construct a developing map of such. 
In the following we assume that all manifolds are smooth and oriented.

\begin{defn} A \emph{tame} manifold is a manifold $M$ diffeomorphic to the interior of a compact manifold $\bar M$. The boundary components of $\bar M$ are called the \emph{ends} of $M$. We allow the number of ends to be zero so that a closed manifold is a tame manifold with no ends.
\end{defn}
To avoid confusing readers familiar with existing terminology, we stress that an end is a boundary component of $\bar M$ and \emph{not} a closed regular neighborhood of a boundary component as many other authors define it. 


Let $M$ be a tame manifold. By the collar neighborhood theorem, we can regard $\bar M$ as a retract of $M$, and we therefore have a canonical identification of $\pi_1(M)$ with $\pi_1(\bar M)$.  
Each of the ends of $M$ defines a subgroup of $\pi_1(M)$, which is well defined up to conjugation. We call these the \emph{peripheral subgroups} of $M$. We neither require that the ends are incompressible nor that the genus is greater than zero, so the peripheral subgroups may be trivial. 
\begin{defn} An element of $G$ is called \emph{parabolic} if it fixes exactly one point in $\partial\bar H^3$. A subgroup of $G$ is called \emph{parabolic} if all its element are parabolic fixing a common point in $\partial\bar H^3$.
\end{defn}
It is easy to see that any parabolic subgroup is conjugate to a subgroup of $P$.
\begin{defn} A representation $\rho\colon \pi_1(M)\to \PSL(2,\C)$ is called \emph{boundary-parabolic} if $\rho$ maps each peripheral subgroup to a parabolic subgroup. An end is called \emph{trivial} with respect to $\rho$ if its corresponding parabolic subgroup is trivial. If $\rho$ is clear from the context we will just call the end trivial.
\end{defn}

\begin{ex}\label{georep} The geometric representation of a hyperbolic $3$--manifold is boundary-parabolic. It is defined up to conjugation. All the ends are non-trivial tori. 
\end{ex}
\begin{ex}\label{galconj}
Recall that any hyperbolic manifold $M$ is isometric to $\H^3/\Gamma$, where $\Gamma$ is a discrete subgroup of $\PSL(2,\C)$. The \emph{trace field} $\Q(\text{tr}\,\Gamma)$ of $M$ is the subfield of $\C$ generated over $\Q$ by the traces of the elements of $\Gamma$. The trace field is a number field, and if $M$ is non-compact, the geometric representation of $M$ is conjugate to a representation into $\PSL(2,\Q(\text{tr}\,\Gamma))$. If $n$ is the degree of $\Q(\text{tr}\,\Gamma)$, there are exactly $n$ embeddings of $\Q(\text{tr}\,\Gamma)$ in $\C$. Composing the geometric representation with the map $\PSL(2,\Q(\text{tr}\,\Gamma))\to \PSL(2,\C)$ induced by one of these embeddings gives a representation which is called a \emph{Galois conjugate} of $\rho$. All Galois conjugates of the geometric representation are boundary-parabolic. We refer to Reid-Maclachlan~\cite{ReidMaclachlan} for more details on trace fields.

\end{ex}
\begin{remark} Suppose $M$ has a single end which is a torus. In this case, the set of conjugation classes of (irreducible) boundary-parabolic representations is often finite. For example, if all components of the $\PSL(2,\C)$--character variety are one dimensional, the characters of boundary-parabolic representations are given by the finite set $I_m^{-1}(4)$, where $m$ denotes a meridian, and $I_\lambda$, for $\lambda\in \pi_1(M)$, is the regular function taking $\rho$ to $(\text{tr}\,\rho(\lambda))^2$. 
In Cooper et al.~\cite{ApolPaper}, the authors prove that all components of the $\SL(2,\C)$--character variety are one dimensional if $M$ contains no closed incompressible surface. The dimension of the $\PSL(2,\C)$--character variety is in general bigger than the $\SL(2,\C)$ analog, but if $M$ is irreducible and if $H_1(M;\Z/2\Z)=\Z/2\Z$ (e.g.~ if $M$ is a knot complement), the dimensions are the same, see Boyer--Zhang~\cite{BoyerZhang}. 
\end{remark}
\begin{remark}\label{twobridge}
If $M$ is a hyperbolic twist knot, it follows from Hoste--Shanahan~\cite{HosteShanahan} that every boundary-parabolic representation, which is not conjugate to a representation in $P$, is a Galois conjugate of the geometric representation.
\end{remark}

\subsection{The developing map of a representation}\label{devsection}

Let $M$ be a tame manifold and let $\hat M$ be the compactification of $M$ obtained by collapsing each boundary component of $\bar M$ to a point. We shall refer to points in $\hat M$ corresponding to the ends as \emph{ideal points} of $M$. Similarly, we let $\hat{\tilde M}$ denote the compactification of the universal cover of $\bar M$ obtained by adding ideal points corresponding to the lifts of the ends of $M$. The covering map extends to a map from $\hat{\tilde M}$ to $\hat M$. In the following we assume that a base point in $M$ and one of its lifts has been fixed once and for all. With the base points fixed we have an action of $\pi_1(M)$ on $\tilde M$ by covering transformations, which extends to an action on $\hat{\tilde M}$. This action is no longer free. The stabilizer of a lift $\tilde e$ of an ideal point $e$ corresponding to an end $E$ is isomorphic to an end subgroup $\pi_1(E)$. Changing the lift $\tilde e$ corresponds to changing the end subgroup by a conjugation. 
\begin{defn}
A \emph{triangulation} of a tame manifold $M$ is an identification of $\hat M$ with a complex obtained by gluing together simplices with simplicial attaching maps. 
We do not require that the triangulation is simplicial but we do require that open simplices embed.
\end{defn}
A triangulation of $M$ always exists. It lifts uniquely to a triangulation of $\tilde M$, and it induces a triangulation of each end of $M$ as the link of the corresponding ideal point. 

\begin{lemma}\label{geodesic} Let $\Delta$ be an $n$--simplex in $\R^n$ with an ordering of the vertices. Given any ideal $n$--simplex $\Delta'\in \H^n$ with a vertex ordering, there exists a unique homeomorphism from $\Delta$ to $\Delta'$ that restricts to an order preserving map of vertices, and takes euclidean straight lines to hyperbolic straight lines.
\end{lemma}
\begin{proof} The existence of such a homeomorphism is obvious if we work in the Klein model of hyperbolic space, where the hyperbolic straight lines and the euclidean straight lines coincide.
The uniqueness follows from the fact that any local homeomorphism between open subsets of $\R^n$ preserving straight lines is affine.
\end{proof} 
We call a homeomorphism as in Lemma~\ref{geodesic} an \emph{ideal homeomorphism}.
\begin{defn}\label{developdef} Let $M$ be a triangulated tame $3$--manifold and let $\rho$ be a boundary-parabolic representation. A \emph{developing map} of $\rho$ is a $\rho$--equivariant map
\begin{equation}D_\rho\colon \hat{\tilde M}\to \bar \H^3\end{equation}
sending all zero-cells to $\partial\bar \H^3$ and satisfying that the composition of $D_\rho$ with the characteristic map of a cell is an ideal homeomorphism onto a non-degenerate ideal simplex. Two developing maps are called \emph{equivalent} if they agree on the ideal points corresponding to lifts of \emph{non-trivial} ends.
\end{defn}
Note that if $D$ is a developing map of $\rho$, then $gD$ is a developing map of $g\rho g^{-1}$.

\begin{thm}\label{devconstr} If the triangulation of $M$ is sufficiently fine, a developing map always exists, and it is unique up to equivalence. A single barycentric subdivision is enough to ensure that any boundary-parabolic representation (including the trivial representation) admits a developing map.
\end{thm}
\begin{proof}
Our construction follows that of Section $8$ of Neumann--Yang~\cite{NeumannYang}. A developing map is uniquely determined by its value on the zero-cells, of which there are three different types to consider: ideal points corresponding to non-trivial ends, ideal points corresponding to trivial ends, and interior zero-cells. 

Let $e$ denote an ideal point of $M$ corresponding to a non-trivial end and let $\tilde e_i$ denote the lifts of $e$. As described above, each of the $\tilde e_i$'s defines a peripheral subgroup. The $\rho$--image of the peripheral subgroup of $\tilde e_i$ is a parabolic subgroup $P_i$ with a unique fixed point $v_i\in \partial \bar \H^3$.  
Define 
\begin{equation}\label{rhoequi}D_\rho(\tilde e_i)=v_i.\end{equation} 
Note that $D_\rho(\alpha \tilde e_i)=\rho(\alpha)v_i$ for every $i$ and every $\alpha\in \pi_1(M)$. 
We define $D_\rho$ on the rest of the zero-skeleton by letting lifts of trivial ends and interior zero-cells map equivariantly to arbitrary points in $\partial\bar \H^3$ requiring that zero-cells in the closure of a $3$--cell map to distinct points. Since we don't allow degenerate simplices, we might need to subdivide the triangulation to get a well defined developing map. This is the case e.g.~if two peripheral subgroups map to the same parabolic subgroup. It is clear that a single barycentric subdivision is enough to ensure non-degeneracy. The uniqueness statement follows from the fact that $\rho$--equivariance forces the image of lifts of non-trivial ideal points to be as in \eqref{rhoequi}. 
\end{proof}

\begin{remark}\label{facepairings} Given a triangulation, the fundamental group is generated by face pairings, and a boundary-parabolic representation is given by an association of an element in $\PSL(2,\C)$ to each such face pairing satisfying the relevant relations. Given this data, the process of developing a boundary-parabolic representation is completely algorithmic and works very fast even for a high number of simplices. 
\end{remark}
Since a developing map of $\rho$ is $\rho$--equivariant, it endows each simplex in the triangulation of $M$ with the shape of an ideal simplex, and it thus allows us to think of $M$ as a space obtained by gluing together ideal simplices. 
If $M$ is a hyperbolic $3$--manifold and $\rho$ is the geometric representation, a developing map provides a degree one ideal triangulation of $M$ in the sense of Neumann--Yang~\cite{NeumannYang}.

In the next section we will define a fundamental class of $\rho$. To obtain this we will need that each of the ideal simplices are decorated, that is, we need to choose euclidean horospheres at each of the ideal vertices. 
\begin{defn}\label{dechorodef} Let $\rho$ be a boundary-parabolic representation of a tame, triangulated $3$--manifold $M$, and let $D$ denote a developing map of $\rho$. 
Let $x\in \hat M$ be a zero-cell. For each lift $\tilde x\in \hat{\tilde M}$ of $x$ let $H(D(\tilde x))$ be a euclidean horosphere based at $D(\tilde x)$. The collection $\{H(D(\tilde x))\}_{\tilde x\in \pi^{-1}(x)}$ of euclidean horospheres is called a \emph{decoration} of $x$ if the following equivariance condition is satisfied:
\begin{equation}\label{equivarhoro}H(D(\alpha \tilde x))=\rho(\alpha)H(D(\tilde x)),\quad \text{for }\alpha \in \pi_1(M),\,\tilde x\in \pi^{-1}(x)\end{equation}
\end{defn} 
\begin{defn}\label{decoration}
Let $M$ and $\rho$ be as above. A developing map of $\rho$ together with a choice of decoration of each zero-cell of $M$ is called a \emph{decoration} of $\rho$. If we have picked a decoration we say that $\rho$ is \emph{decorated}. We consider two decorations to be equivalent if they agree on non-trivial ideal points.
\end{defn}
Note that a decorated representation endows each simplex of $\hat M$ with the structure of a decorated ideal simplex determined up to congruence.

\section{The fundamental class of a representation}\label{fundclasssec}
In this section we show that the notion of decoration from Definition \ref{decoration} extends to the more general setup of $(G,H)$--representations. We show that a decorated $(G,H)$--representation determines a fundamental class in $H_*(G,H)$, and we describe a particularly simple way of constructing this class in the special case of boundary-parabolic representations. The general theory seems interesting in itself, and we develop it in detail. 
All manifolds are assumed to be smooth and oriented, but not necessarily of dimension $3$.

\begin{defn} Let $M$ be a tame manifold and let $H$ be a subgroup of $G$, where $G$ is any (discrete) group. A representation $\rho\colon \pi_1(M)\to G$ is called a $(G,H)$--\emph{representation} if $\rho$ sends peripheral subgroups to conjugates of subgroups of $H$.
\end{defn}  

\subsection{Definition of the fundamental class}
We start with the case where $M$ is closed and $H$ is the trivial subgroup. In this case a $(G,H)$--representation is just a representation. Recall that conjugation classes of representations of $\pi_1(M)$ are in one-one correspondence with homotopy classes of \emph{classifying maps} $M\to BG$.
\begin{defn} Let $\rho\colon \pi_1(M)\to G$ be a representation and let $f$ denote its classifying map. The \emph{fundamental class} of $\rho$ is the class $f_*[M]$, where $[M]$ is the fundamental class of $M$. 
\end{defn} 

Suppose that $M$ is triangulated. As mentioned in Section
\ref{relative} a generator of $B_*(G)$ can be regarded as a simplex
with a $G$--labeling of edges. We can therefore produce cycles in
$B_*(G)$ by labeling the edges of $M$ in an appropriate fashion. The definition below is
taken from Neumann~\cite{Neumann}. 

\begin{defn} Let $M$ be a triangulated manifold. Let $S_q(M)$ be the set of oriented $q$--cells of $M$. A $G$--\emph{cocycle} on $M$ is a map $\sigma\colon S_1(M)\to G$ satisfying the properties:
\begin{enumerate}[(i)]
\item $\sigma\langle v_0,v_2\rangle=\sigma\langle v_0,v_1\rangle\sigma\langle v_1,v_2\rangle\quad \text{for}\quad\langle v_0,v_1,v_2\rangle\in S_2(M)$.
\item $\sigma\langle v_1,v_0\rangle = \sigma\langle v_0,v_1\rangle^{-1}$.
\end{enumerate}
If $\tau\colon S_0(M)\to G$ is a $0$--cochain, then its \emph{coboundary action} on $G$--cocycles is to replace $\sigma$ by 
\begin{equation}\label{cobact}\langle v_0,v_1\rangle\mapsto \tau(v_0)^{-1}\sigma\langle v_0,v_1\rangle \tau(v_1).\end{equation}
\end{defn}

A $G$--cocycle $\sigma$ gives rise to a representation $\rho\colon
\pi_1(M)\to G$ which is well defined once we have chosen a zero-cell
as a base point. We say that $\sigma$ \emph{represents} $\rho$. Given a representation we can always find a cocycle representing it, e.g.~by defining the cocycle to be the identity on edges of a maximal tree. A representing cocycle is unique up to the action by coboundaries. From this the proposition below easily follows. It does not require that $M$ be closed.
\begin{prop}\label{cocyclehomo}  There is a one-one correspondence
  between $G$--cocyc\-les up to the action by coboundaries, and homotopy
  classes of classify\-ing maps $M\to BG$.
\end{prop}
To obtain a cycle in $B_*(G)$ from a $G$--cocycle, we need that each
simplex in the triangulation of $M$ has a vertex ordering which is
respected by the face identifications.
\begin{defn}\label{ordering}An \emph{ordering} of a triangulated tame manifold $M$ is an ordering of the vertices of each simplex satisfying that the orientation of edges induced by the ordering agrees under the identification of faces. Having fixed an ordering we say that $M$ is \emph{ordered}.
\end{defn}
\begin{remark}Not every triangulation has an ordering, but after performing a single barycentric subdivision, we have a natural ordering by codimension. Namely, the $i$--th vertex of a simplex is the unique vertex lying in a face of codimension $i$ in the original simplex.\end{remark}

\begin{prop}\label{expgen}
Let $M$ be a closed, ordered, triangulated $n$--manifold with a representation $\rho\colon \pi_1(M)\to G$. Let $\sigma$ be a
$G$--cocycle representing $\rho$ and let $\Delta_i$ be the simplices of
$M$ endowed with the $G$--labeling of oriented edges induced by $\sigma$. 
Let $\epsilon_i$ be a sign indicating whether or not the orientation of $\Delta_i$ 
induced by the ordering agrees with the orientation it inherits from $M$.
The cycle 
\begin{equation}\label{classlabel}
\sum \epsilon_i \Delta_i\in B_n(G)\end{equation}
represents the fundamental class. 
\end{prop}
\begin{proof}
The proof is an application of the Milnor construction of $BG$, which we recall below. 
Let $\Delta^n=\{(t_0,\dots,t_n)\in \R^{n+1}\mid \sum_i t_i=1\}$ be the standard simplex and let $\partial_i$ denote the map $\Delta^{n-1}\to \Delta^n$ inserting a zero on the $i$--th coordinate. 
We have
\begin{equation}\label{BGrelation}
BG=\left(\bigsqcup_{n=1}^\infty \Delta^n\times G^n\right) \Big / \sim
\end{equation}
where the relation is generated by $(\partial_i t,x)\sim (t,d_i x)$, with 
\[d_i(g_1,\dots g_n)=\begin{cases}(g_2,\dots ,g_n)&\text{for }i=0\\(g_1,\dots,g_ig_{i+1},\dots,g_n)&\text{for } 0<i<n\\
(g_1,\dots,g_{n-1})&\text{for }i=n\end{cases}.\]
Note that the set $G^n$ parametrizes the set of $G$--cocycles on $\Delta^n$. Namely, a tuple $(g_1,\dots,g_n)$ corresponds to the unique cocycle sending the edge $[e_{i-1},e_i]$ to $g_i$, where $e_i$ is the $i$--th standard basis vector. This means that an ordered simplex with a $G$--cocycle is naturally equipped with a map to $BG$. The cocycle $\sigma$ induces a $G$--cocycle of each simplex of $M$, and since $M$ is ordered, the maps to $BG$ respect the face pairings, and thus induce a map from $M$ to $BG$. By construction, the induced map on $\pi_1$ is $\rho$, so it is a classifying map. Using the canonical isomorphism between the cellular complex of $BG$ with the complex $B_*(G)$, it follows that the fundamental class has the given representation.
\end{proof}

We now return to the general case of a $(G,H)$--representation of a tame manifold. 
Recall that a triangulation of $M$ induces a triangulation of
$\partial \bar M$. It also induces a cell decomposition of $\bar M$ using
hybrids of truncated simplices and normal simplices. A simplex of $M$ with vertices consisting entirely of ideal points gives rise to a truncated simplex, and a simplex of $M$ consisting entirely of interior points gives rise to a normal simplex.
The notion of a $G$--cocycle
extends to such cell decompositions in the obvious way.

\begin{defn}
A $(G,H)$--cocycle is a $G$--cocycle on $\bar M$ sending edges of
$\partial \bar M$ to $H$. A $0$--cochain sending vertices of $\partial
\bar M$ to $H$ acts on a $(G,H)$--cocycle as in \eqref{cobact}.
\end{defn}

\begin{prop}\label{reldiag}
There is a one-one correspondence between $(G,H)$--co\-cycles up to the action by coboundaries, and homotopy commuting diagrams
\begin{equation}\label{boundarydiag}\cxymatrix{{\partial \bar M\ar[d]\ar[r] &BH\ar[d]\\\bar M\ar[r]&BG.}}\end{equation}
\end{prop}
\begin{proof} Using the Milnor construction of $BG$ we see that a $(G,H)$--cocycle induces a diagram as above. Now suppose we have a diagram as above. Let $f$ denote the map from $\bar M$ to $BG$, and let $\sigma$ be a $G$--cocycle representing $f$. Since the restriction of $f$ to $\partial \bar M$ is homotopy equivalent to a map into $BH$, Proposition~\ref{cocyclehomo} implies that the restriction of $\sigma$ to $\partial \bar M$ can be modified by coboundaries to send edges of $\partial \bar M$ to $H$. The $(G,H)$--cocycle thus produced is unique up to the action by coboundaries, and the induced map is homotopy equivalent to $f$.
\end{proof}
We can now define the fundamental class of a $(G,H)$--cocycle as
the image of the fundamental class $[\bar M,\partial \bar M]$ under
the corresponding map \eqref{boundarydiag}.

Let $\bar M'$ be the manifold obtained from $\bar M$ by removing small
disjoint open balls around each interior zero-cell. 
Note that the triangulation of $M$ induces a cell decomposition of $\bar M'$
consisting entirely of truncated simplices. If $M$ has dimension $n\geq 3$, it follows from Proposition~\ref{reldiag} that there is a one-one correspondence between cocycles on $\bar M$ and cocycles on $\bar M'$ (modulo coboundaries), and that corresponding cocycles induces the same fundamental class.
We have the following generalization of Proposition~\ref{expgen}. 

\begin{prop}\label{truncexplicit}
Let $M$ be an ordered, triangulated manifold of dimension $n\geq 3$, and let $\tau$ be a $(G,H)$--cocycle on $\bar M'$. 
The fundamental class is represented in $\bar B_n(G,H)$ by the cycle
\[\sum \epsilon_i\bar\Delta_i,\]
where $\bar \Delta_i$ are the truncated simplices in the triangulation of $\bar M'$ with edge labelings given by $\tau$, and $\epsilon_i$ is a sign indicating whether or not the orientation of $\bar \Delta_i$ given by the ordering agrees with the orientation induced from $\bar M'$.
\end{prop} 

\begin{proof}
The proof is similar to the proof of Proposition \ref{expgen}. We give a sketch and leave the details to the reader.
Let $T^n$ be a set parametrizing the set of $(G,H)$--cocycles of a truncated $n$--simplex. Consider the space 
\[B(G,H)=\left(\bigsqcup_{n=1}^\infty \bar\Delta^n\times T^n\right) \Big / \sim\]
where the relation is the truncated analog of the relation in \eqref{BGrelation}. The map from $G^n$ to $T^n$ obtained by labeling all short edges by $1$ induces a map from $BG$ to $B(G,H)$, whose restriction to $BH$ can be seen to be null homotopic. By the universal property of the homotopy cofiber, $B(G,H)$ is a model for the cofiber of $BH\to BG$.
The homology of $B(G,H)$ is equal to the homology of $\bar B_*(G,H)$, and the fundamental class is easily seen to have a representation as in the proposition.
\end{proof}

\subsection{The $(G,H)$ cocycle}
We will now describe how to
associate a $(G,H)$--cocycle to a $(G,H)$--representation $\rho$. By Proposition~\ref{reldiag}
this defines a fundamental class. The cocycle, and therefore
also the fundamental class, will depend on a choice of decoration which we define below.

Let $e_i$ denote the ideal points corresponding to the non-trivial ends $E_i$ of $M$, and choose lifts $\tilde e_i\in \hat{\tilde M}$ of $e_i$. As in Section \ref{devsection} this defines peripheral subgroups $\pi_1(E_i)$. Let $H_i$ denote the image of $\pi_1(E_i)$. Note that replacing $\tilde e_i$ by $\alpha \tilde e_i$ replaces $H_i$ by its conjugate $\rho(\alpha)H_i\rho(\alpha)^{-1}$.
Pick a $G$--cocycle of $\bar M$ that represents $\rho$ and sends edges of $E_i$ to $H_i$. 
To see that such a cocycle exists, we can construct it as follows:
pick an $H_i$--cocycle $\sigma_i$ on $E_i$ representing the restriction of $\rho$ to $\pi_1(E_i)$. Define a cocycle on $\bar M$, by letting its restriction to $E_i$ be $\sigma_i$ and letting it be $1$ on edges of a maximal tree in $\bar M\backslash\cup_i E_i$ containing edge paths from the base point in $M$ to $E_i$. This uniquely specifies the value on all edges of $\bar M$.
Up to the action by coboundaries sending vertices of $E_i$ to $H_i$, this cocycle is unique once the choices of lifts $\tilde e_i$ have been fixed. 
Pick elements $g_i$ satisfying $g_i^{-1}H_ig_i\subset H$. We call these \emph{conjugation elements}.
Modifying the above $G$--cocycle by the coboundary of the $0$--cochain
\[\tau(v)=\begin{cases}g_i \text{ if }v\in E_i\\1 \text{ otherwise}\end{cases}\]
gives us a $(G,H)$--cocycle, which up to the action by coboundaries, only depends on the choices of conjugation elements. We will refer to it as the $(G,H)$--cocycle \emph{associated} to $\rho$. Note that multiplying a conjugation element from the right by an element in $H$, changes the $(G,H)$--cocycle by a coboundary. In the following we will thus regard the conjugation elements as left $H$--cosets.

The conjugation elements, and therefore also the associated $(G,H)$--cocycle, depend on the choices of lifts. To indicate this dependence we will now denote them $g_i(\tilde e_i)$. To make the associated $(G,H)$--cocycle independent of the choices of lifts, the conjugation elements have to be chosen in an equivariant fashion.

\begin{defn}\label{decoconjel} Let $\rho$ be a $(G,H)$--representation.
A set of conjugation elements $g_i(\alpha \tilde e_i)$, $\alpha\in \pi_1(M)$, satisfying the equivariance condition
\begin{equation}\label{equivarconj}g_i(\alpha \tilde e_i)=\rho(\alpha)g_i(\tilde e_i),\quad \alpha\in \pi_1(M).
\end{equation}
 is called a \emph{decoration} of $\rho$. 
\end{defn}

\begin{remark} Note that decorations are parametrized by the group $\big(N_G(H)/H\big)^n$, where $n$ is the number of non-trivial ends and $N_G(H)$ is the normalizer of $H$ in $G$.
\end{remark}
Given a decoration, the associated $(G,H)$--cocycle is well defined, and unique up to the action by coboundaries. 
We have thus proved:
\begin{thm}\label{fundclassexist} A decoration of a $(G,H)$--representation determines a fundamental class in $H_*(G,H)$.
\end{thm}

\subsection{An explicit construction of the fundamental class}
We now specialize to the case of boundary-parabolic representations of tame $3$--manifolds, i.e.~the case with $G=\PSL(2,\C)$, and $H=P$. 
In this case there is a particularly simple way of constructing the $(G,P)$--cocycle. 


\begin{lemma}\label{agree} For a boundary-parabolic representation of a tame $3$--manifold, there is a natural one-one correspondence between decorations by euclidean horospheres and decorations by conjugation elements.
\end{lemma}
\begin{proof}   
By Remark \ref{Hinf}, a euclidean horosphere $H(v)$ at $v\in \partial \bar \H^3$ corresponds to a left coset $gP$, where $g$ takes $H(\infty)$ to $H(v)$. Hence, we only need to check that the two notions of equivariance as defined in Definition \ref{dechorodef} and Definition \ref{decoconjel} agree. We leave this to the reader. 
\end{proof} 

Let $\rho$ be a decorated boundary-parabolic representation of a tame $3$--manifold $M$. We will assume that $M$ is ordered. This is no restriction since we can always obtain an ordering by performing a barycentric subdivision.
Recall that $\rho$ endows each of the $3$--cells of $M$ with the shape of a decorated ideal simplex. 
We can thus think of $M$ as a collection $\{\Delta_i\}$ of decorated
ideal tetrahedra together with a set of face pairings. By Theorem~\ref{trunconetoone} each $\Delta_i$ corresponds to a generator $\bar\Delta_i$ of $\bar B_3(G,P)$, which is a truncated simplex together
with a labeling of its oriented edges. 
The face pairings of $\Delta_i$ induce face pairings of the
corresponding truncated simplices $\bar \Delta_i$. Note that the complex obtained by gluing these together is homeomorphic to $\bar M'$, the manifold obtained from $\bar M$ by removing disjoint open balls around each interior zero-cell.  Since both the decorations and the orderings respect the face pairings, we see that the face pairings of the truncated
simplices respect the labelings of short edges. Since the labelings
of long edges are determined uniquely by the labelings of the short
edges, the face pairings respect the labelings of long edges as well.
This means that the labelings form a $(G,P)$--cocycle $\sigma$ of $\bar M'$.

\begin{thm}\label{representsfund} The cocycle $\sigma$ is the $(G,P)$--cocycle associated to $\rho$.
\end{thm}
\begin{proof}
Up to multiplication by elements of $G$, a developing map can be uniquely reconstructed from the ideal simplex shapes and the gluing pattern. Since a representation is determined by the equivalence class of its developing map, it follows from Theorem \ref{trunconetoone} that $\sigma$ represents $\rho$. The fact that $\sigma$ is the $(G,P)$--cocycle associated to the correct decoration of $\rho$ is an easy consequence of the observation in Remark \ref{newviewpoint}.
\end{proof}

\begin{cor}
The cycle \begin{equation}\label{fundclassrep}
\sum \epsilon_i\bar \Delta_i\in\bar B_3(G,P)
\end{equation} represents the fundamental class. As always, $\epsilon_i$ is a sign which is positive if and only if the orientation given by the ordering agrees with the induced orientation.
\end{cor}

\begin{remark}
If $M$ is not ordered the above construction still works, but we need the ordering to represent the fundamental class in the truncated complex. Performing a barycentric subdivision produces new non-ideal zero cells, and a decoration can be obtained by equivariantly ``shooting off'' these new zero cells to random points in $\partial \bar \H^3$ and picking random decorations at these points satisfying the equivariance condition \eqref{equivarhoro}. As before, this produces a $(G,P)$--cocycle $\sigma$. The ``shooting off'' process changes the geometry, but the fundamental class, which only depends on the decoration at ideal points, is still given by $\sigma$.
\end{remark}

\begin{ex}\label{knotex} We illustrate the above by constructing the fundamental class of the geometric representation of the $5_2$ knot complement. This manifold has an ordered triangulation as shown in Figure~\ref{knotcomp}. For each of the simplices the orientation determined by the ordering agrees with its inherited orientation, i.e., the signs $\epsilon_i$ are all positive. 
\begin{figure}[ht]
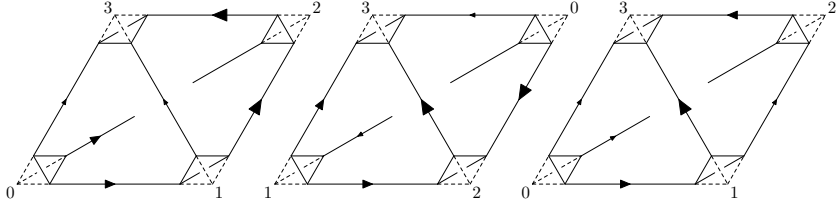
\label{gluingpattern}
\centering
\begin{minipage}[c]{0.3\textwidth}
\centering 
\includegraphics[width=4.2cm]{ChernSimonsFigure.5}
\end{minipage}
\negthickspace\negthickspace\negthickspace\negthickspace\negthickspace
\begin{minipage}[c]{0.3\textwidth}
\centering 
\includegraphics[width=4.2cm]{ChernSimonsFigure.6}
\end{minipage}
\negthickspace\negthickspace\negthickspace\negthickspace\negthickspace
\begin{minipage}[c]{0.3\textwidth}
\centering 
\includegraphics[width=4.2cm]{ChernSimonsFigure.7}
\end{minipage}
\caption[Gluing pattern for the $5_2$ knot complement]{Gluing pattern for the $5_2$ knot complement.}\label{knotcomp}
\end{figure}
As described e.g.~in Neumann--Zagier~\cite{NeumannZagier}, the ideal simplex shapes obtained by developing the geometric representation can be found by solving a set of polynomial equations called the \emph{gluing equations}. In this case there are five gluing equations, one for each of the edges and two for the cusp. The gluing equation for an edge states that the product of the cross-ratio parameters associated to it is $1$. Letting $u,v$ and $w$ denote the cross-ratios, we can read off the gluing equations for the edges from Figure \ref{knotcomp}. We obtain:
\[u'u''vv'w'^2w''=1,\quad uu''v'v''w^2=1,\quad uu'vv''w''=1.\]
Figure \ref{cuspdev} below shows a developing of the cusp. 
\begin{figure}[ht]
\centering
\includegraphics[width=6.5cm]{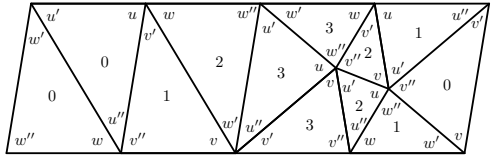}
\caption[Developing image of the cusp]{Developing image of the cusp.}\label{cuspdev}
\end{figure} 
Since the cusp is a torus, the peripheral subgroup has two generators, and the requirement that these both map to parabolics yields two equations, which can be read off from Figure~\ref{cuspdev}:
\[w^{-1}u'=1,\quad w'u'uv'ww''u'w'wv'uu''v'=1.\]
The ideal triangulation is given by the unique solution of the above five equations satisfying that $u$, $v$ and $w$ all have positive imaginary part.
After a little algebraic manipulation we obtain that the solutions are given by $u=x^2$, $v=x^2-x+1$ and $w=v=x^2-x+1$, where $x$ is a root of $x^3-x^2+1$. The geometric solution corresponds to the root with positive imaginary part. The two other solutions correspond to Galois conjugates.

A decoration allows us to view the configuration in Figure~\ref{cuspdev} as a configuration in $\C$, and we can thus assign a complex number to each oriented edge as defined by the vector going from the start point to the end point. In this example we choose the decoration such that the lower left edge has vertices at $0$ and $1$. Using Figure~\ref{Ctriangle} we see that the labelings of all edges are given by products of cross-ratio parameters, e.g.~the shortest edge (oriented downwards) is labeled by $w(v'')^{-1}=u^{-1}v$. Each of the triangles in Figure~\ref{cuspdev} corresponds to a cut-off triangle as indicated by the numbers.
This gives us labelings of the short edges in Figure~\ref{knotcomp} which obviously satisfy \eqref{alphas}, and by 
Lemma~\ref{unique}, the labelings of the long edges are uniquely determined by \eqref{squareroot}.
We obtain that the edges marked with the small arrow, the bigger arrow and the biggest arrow, respectively, are labeled by the matrices 
\begin{equation*}\Mat{0}{-a^{-1}}{a}{0},\quad \Mat{0}{-b^{-1}}{b}{0},\quad \Mat{0}{-c^{-1}}{c}{0},\end{equation*} 
with $a$, $b$ and $c$ satisfying
\begin{equation}\label{matrices}
a^2=(w'')^{-1}=1-x,\quad b^2=1,\quad c^2=(ww'^{-1})^{-1}=x^2-x-1.\end{equation} 

This finishes the construction of the $(G,P)$--cocycle associated to the geometric representation, and the fundamental class is now given by \eqref{fundclassrep}.
\end{ex}
Note that all we need to produce the fundamental class is the gluing pattern of the triangulation together with the developing image of each of the ends.

\begin{remark} To see the dependence of the decoration explicitly, consider the natural map $H_3(\PSL(2,\C),P)\to H_2(P)$. An explicit formula for this map is given by the map taking a truncated simplex to the sum of the four triangles determined by the small edges. Changing the decoration corresponds to changing the labelings of short edges by an element in $\C^*$, and it thus follows that different decorations yield fundamental classes having different images in $H_2(P)$. In the general case, the fundamental class depends on the decoration whenever $N_G(H)/H$ acts non-trivially on the homology of $H$. 
\end{remark}

\subsection{Other pairs}

We briefly discuss how the construction of the fundamental class works for other pairs of groups.
\subsubsection{Boundary-loxodromic representations}
Suppose $M$ is a tame $3$--manifold where all the ends of $M$ are tori. A $G$--representation of $\pi_1(M)$ taking peripheral subgroups to subgroups fixing a unique geodesic, is a $(G,T)$--representation, where $T$ is the subgroup of diagonal elements. Since $N_G(T)/T=\Z/2\Z$, there are $2^n$ different decorations, where $n$ is the number of ends. By Theorem \ref{fundclassexist} each of these gives rise to a fundamental class in $H_3(G,T)$. 
One can prove that $H_3(G,T)=\Pre (\C)$, so the fundamental classes of boundary-loxodromic representations contain no more information than the cross-ratios.
\subsubsection{Higher dimensions}
Consider the pair $(\SL(n,\C),P)$, where $P$ is the subgroup of $\SL(n,\C)$ consisting of upper triangular matrices with $1$ on the diagonal.
Let $M^k$ be a tame $k$--manifold, $k\geq 3$. 
By Theorem \ref{fundclassexist}, a decorated $(\SL(n,\C),P)$--representation of $\pi_1(M)$ gives rise to a fundamental class in $H_k(\SL(n,\C),P)$. 
We wish to find an explicit representative.
The complex $\bar B_*(\PSL(2,\C),P)$ from Section~\ref{compoftrunc} generalizes in an obvious way to a complex $\bar B_*(\SL(n,\C),P)$, based on labelings of truncated simplices, where long edges are labeled with counter diagonal elements and short edges by elements in $P$.
The explicit construction is via a generalization of Lemma~\ref{uniquerep}. We state it for $n=3$. The general case is similar.

For $g\in \SL(3,\C)$, let $M_{ij}(g)$ be the minor obtained from $g$ by removing the $i$th row and the $j$th column. As usual, we let $g_{ij}$ denote the $ij$th entry of $g$.
\begin{lemma}\label{higherdim} Let $gP$, $hP$ be such that 
\[(g^{-1}h)_{31}\neq 0\text{  and  }\det (M_{13}(g^{-1}h))\neq 0.\] Note that this is independent of the representatives.
There are unique coset-represen\-tatives $gu$ and $hv$ such that $(gu)^{-1}hv$ is counter diagonal. Explicitly, if 
\[g^{-1}h=\begin{pmatrix}a&b&c\\d&e&f\\g&h&i\end{pmatrix},\,\, u=\begin{pmatrix}1&x&y\\0&1&z\\0&0&1\end{pmatrix},\,\, 
v=\begin{pmatrix}1&r&s\\0&1&t\\0&0&1\end{pmatrix},\,\, \delta=\begin{pmatrix}0&0&\gamma\\0&\beta&0\\\alpha&0&0\end{pmatrix}\]  
then $(gu)^{-1}hv=\delta$ if and only if 

\begin{equation}
\begin{gathered}
r=\frac{-h}{g},\quad y=\frac{a}{g},\quad  z=\frac{d}{g},\\
x=\frac{ah-bg}{dh-eg},\quad  s=\frac{ei-fh}{dh-eg},\quad  t=\frac{fg-di}{dh-eg},\\
\alpha=g,\quad \beta=\frac{-(dh-eg)}{g},\quad \gamma=\frac{1}{dh-eg}.
\end{gathered}
\end{equation}
\end{lemma}
A decoration of $\rho$ associates left $P$--cosets to each vertex of a fundamental domain of $M$ in $\tilde M$. If we view $\pi_1(M)$ as being generated by face pairings of the fundamental domain, this process is completely explicit. Given the $P$--cosets, a representative for the fundamental class in $\bar B_k(\SL(3,\C),P)$ can be constructed explicitly using Lemma~\ref{higherdim}. We leave the details to the reader.

\section{The complex volume}\label{compvolsec}

We return to the case $G=\PSL(2,\C)$. Even though the fundamental class depends on the choice of decoration, it turns out that its image in $\widehat \B(\C)$ under the map $\Psi\colon H_3(G,P)\to \widehat\B(\C)$ does not. 


Recall that a closed curve in a triangulated complex is called \emph{normal} if it does not intersect the one-skeleton, and intersects every two-cell that it meets transversely. Let $K$ be a complex obtained by gluing together ideal simplices.

 The definition below is a slight generalization of Neumann~\cite[Definition~$4.4$]{Neumann}.
\begin{defn} A \emph{weak flattening} of $K$ is a flattening of each ideal simplex of $K$ such that the total log-parameter around each edge is zero. If the total log-parameter along any normal curve in the star of each zero-cell is zero, it is called a \emph{semi-strong flattening}. The log-parameters must be summed according to the sign conventions of Neumann~\cite[Definition~$4.3$]{Neumann}. 
\end{defn}

\begin{defn} The \emph{parity} of an edge $E$ in a simplex of $K$ is defined by $r \text{ mod }2$, where $w(E)=\Log(z(E))+r\pi i$, and $w(E)$ and $z(E)$ are the log-parameter and cross-ratio parameter of $E$.
\end{defn}

\begin{defn}
A semi-strong flattening satisfying that the total parity parameter is zero along any normal curve in $K$, is called a \emph{strong} flattening.
\end{defn}

The theorem below summarizes the main results of Neumann \cite{Neumann}.

\begin{thm}[\textbf{Neumann}]\label{Lhat} There is a canonical isomorphism $\lambda\colon H_3(G)\cong \widehat \B(\C)$ satisfying that $\hat c_2=\widehat L\circ \lambda$. Furthermore, if $M$ is a cusped hyperbolic manifold $M$ with an ideal triangulation, any strong flattening of $M$ determines a fundamental class $\alpha\in \widehat \B(\C)$ satisfying
\begin{equation}
\widehat L(\alpha)=i(\Vol(M)+i\CS(M)). 
\end{equation}
\end{thm}
Recall the map $\Psi\colon H_3(G,P)\to \widehat\B(\C)$. By Theorem~\ref{trunconetoone} a class in $H_3(G,P)$ can be represented by a collection of decorated ideal simplices, and in this picture, we can view $\Psi$ as a way of endowing each of these simplices with a flattening. 

\begin{figure}[ht]
\centering
\includegraphics[width=6cm]{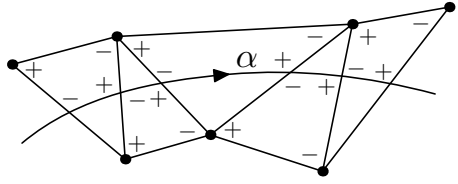}
\caption{A normal curve $\alpha$ in the star of a zero-cell.}\label{curve}
\end{figure}

\begin{thm}\label{strongflat} 
Let $K$ be a complex of decorated ideal simplices representing a homology class in $H_3(G,P)$. The map $\Psi$ endows $K$ with a semi-strong flattening.
\end{thm} 
\begin{proof} Consider a normal curve $\alpha$ in the star of a zero-cell. 
Figure~\ref{curve} shows $\alpha$ as viewed from the zero-cell. Each triangle in the figure corresponds to an ideal simplex. A vertex $v$ of a triangle corresponds to an edge $e$ of the simplex, and the side opposite $v$ corresponds to the edge opposite $e$. Note that some of the triangles might be ``folded back'' or ``flat'' even though this is not indicated in the figure. 
Whenever $\alpha$ passes through a simplex $\Delta$ it picks up a log-parameter of an edge. By \eqref{logpar}, the log-parameter of an edge $e$ is a signed sum of the $\Log(c)$-parameters associated to the four edges of $\Delta$ that are neither $e$ nor the edge opposite $e$. The signs are indicated in the figure. Note that the $\Log(c)$--parameters are associated to the one-cells of $K$, that is, edges that are identified in $K$ have the same $\Log(c)$-parameter. A side of a triangle thus has three $\Log(c)$-parameters associated to it (one for each vertex and one for the side). From Figure~\ref{curve} we see that the total log-parameter along the curve $\alpha$ is a signed sum of $\Log(c)$-parameters associated only to the sides where $\alpha$ enters and exits. The other $\Log(c)$-parameters cancel out. From this we conclude that if $\alpha$ is a closed curve, the total log-parameter along $\alpha$ is zero. 
\end{proof}

\begin{remark} It is not true in general that $\Psi$ provides a strong flattening. 
\end{remark}

\begin{remark}\label{parityhomo}
Neumann shows that parity along normal curves can be viewed as a cohomology class in $H^1(M;\Z/2\Z)$, and by \cite[Corollary 5.4]{Neumann}, the parity of a semi-strong flattening is an element in $\Ker\big(H^1(M;\Z/2\Z)\to H^1(\partial \bar M;\Z/2\Z)\big)$.
\end{remark}

We wish to prove that the image of the fundamental class in $\widehat \B(\C)$ is independent of the decoration. To obtain this, we need to recall some combinatorics of $3$--cycles used in Neumann~\cite{Neumann}. This first appeared in Neumann~\cite{NeumannComb}.

Let $M$ be a tame $3$--manifold with an ordered triangulation. For each $3$--simplex $\Delta$ of $M$, we denote the six edges of $\Delta$ by $e_i$ in the following way: $e_0$ is the edge between vertex $0$ and vertex $1$ of $\Delta$, $e_1$ is the edge between vertices $1$ and $2$, $e_2$ is the edge between vertices $0$ and $2$. For $i=3,4,5$, $e_i$ is the edge opposite $e_{i-3}$. We associate to $\Delta$ a $\Z$--module $J_\Delta$ generated by the $e_i$'s and subject to the relations
\begin{gather}\label{Jdelta}
\begin{aligned} e_i-e_{i+3}=0 \text{ for }i=0,1,2.\\
e_0+e_1+e_2=0.\end{aligned}
\end{gather}
Let $J$ be the direct sum $\bigoplus J_\Delta$, summed over all the $3$--simplices of $M$. For $i=0,1$ let $C_i$ be the free $\Z$--module on the unoriented $i$--cells of $M$.
There is a chain complex
\begin{equation*}
\xymatrix{0\ar[r]&C_0\ar[r]^-\alpha&C_1\ar[r]^-\beta&J\ar[r]^{\beta^*}&C_1\ar[r]^{\alpha^*}&C_0\ar[r]&0}
\end{equation*} 
The map $\alpha$ takes a vertex to the sum of the incident $1$--cells, and the map $\beta$ is defined such that its $J_\Delta$--component takes a $1$--cell $E$ of $M$ to the sum of the edges $e_i$ of $\Delta$ which are identified with $E$ in $M$. We will not need the definition of the other maps.

Recall that a boundary-parabolic representation of $M$ allows us to regard $M$ as a complex of ideal simplices. Suppose we have picked flattenings $(w^i_0,w^i_1,w^i_2)$ of each ideal simplex $\Delta_i$. Consider the element $w\in J\otimes\C$ whose $\Delta_i$--component is 
\begin{equation}\label{flatj}\epsilon_i(w_1^ie_0-w_0^ie_1)\in J_{\Delta_i}\otimes \C.
\end{equation}
Here $\epsilon_i$ is a sign indicating whether or not the vertex-ordering of $\Delta_i$ agrees with the orientation inherited from $M$. By an abuse of notation we will for $j\in J\otimes \C$ write $\epsilon_i j$ for the element in $J\otimes\C$ whose $\Delta_i$--component is $\epsilon_i$ times the $\Delta_i$--component of~$j$.
The following lemma collects results shown in Neumann~\cite{Neumann}.
\begin{lemma}\label{betastar} The flattenings constitute a weak flattening of $M$ if and only if $\beta^*(w)=0$. In this case, the corresponding element in $\widehat \Pre(\C)$ lies in $\widehat \B(\C)$. Furthermore, this element is invariant under replacing $w$ by $w-\pi i\epsilon_i\beta(x)$, where $x$ is an element in $C_1\otimes \Z$.
\end{lemma}
Recall that a decoration associates $\Log(c)$--parameters to each $1$--cell of $M$. It thus gives rise to an element $l$ in $C_1\otimes \C$. The corresponding flattenings given by \eqref{logpar} give rise to an element $w$ in $J\otimes \C$ as above.
\begin{lemma}\label{betaeqw} The elements $l\in C_1\otimes \C$ and $w\in J\otimes \C$ as above are related by $\beta(l)=-\epsilon_i w$. 
\end{lemma}
\begin{proof} Let $l_i$ be the $e_i$--coefficient of $\beta(l)$ in $J_\Delta\otimes \C$. By definition of the map $\beta$, $l_i$ is the $\Log(c)$--parameter of the $1$--cell corresponding to $e_i$. Using \eqref{Jdelta}, we 
have 
\begin{equation*}\beta(l)=l_0e_0+\dots +l_5e_5=(l_0+l_3-l_2-l_5)e_0+(l_1+l_4-l_2-l_5)e_1,
\end{equation*}
which by \eqref{logpar} equals the $\Delta$--component of $-\epsilon_i w$.
\end{proof}

\begin{thm}\label{choiceoflog} The image of the fundamental class in $\widehat\B(\C)$ is independent of the choice of decoration. 
\end{thm}
\begin{proof} For $i=1,2$, let $l_i\in C_1\otimes \C$ and $w_i\in J\otimes \C$ be defined by two different decorations.
By Lemma \ref{betaeqw} we have $\beta(l_2-l_1)=\epsilon_i(w_1-w_2)$. Since $l_2-l_1=\alpha(x)+\pi i y$ for some $y\in C_1\otimes \Z$ and some $x\in C_0\otimes \C$, Lemma \ref{betastar} implies that the elements in $\widehat\B(\C)$ are the same.
\end{proof}

\begin{remark}\label{choice}
By a similar argument, the map $\Psi$ is independent of the branch of logarithm used in \eqref{logpar} to define the flattenings.\end{remark}

Recall that $H_3(G)$ is canonically isomorphic to $\widehat \B(\C)$. In \cite[Proposition 14.3]{Neumann} Neumann shows that the long exact sequence for the pair $(BP,BG)$ gives rise to a split exact sequence
\begin{equation}\label{splitting}\xymatrix{0\ar[r]& H_3(G)\ar[r]^-{i_*}& H_3(G,P)\ar[r]^-{\partial_*}&H_2(P)\ar[r]&0.}\end{equation}

\begin{prop}\label{splittingprop} 
Identifying $H_3(G)$ with $\widehat\B(\C)$, the map $\Psi$ defines a splitting of the sequence \eqref{splitting}. 
\end{prop}
\begin{proof} 
By Neumann~\cite[Lemma $11.3$]{Neumann}, the failure of the parity condition effects the element in $\widehat \B(\C)$ at most by the unique element in $\widehat \B(\C)$ of order $2$. This means that the homomorphism $\Psi\circ i_*-\id$ has image of order at most $2$. Since $H_3(G)$ is divisible, it can have no non-trivial finite quotient. Hence, $\Psi\circ i$ equals the identity.
\end{proof}

\begin{remark}\label{semi} It follows from Proposition \ref{splittingprop} that the semi-strong flattening given by $\Psi$ gives rise to the same element in $\widehat \B(\C)$ as a strong flattening. This is not true for an arbitrary semi-strong flattening. The key point is that $\Psi$ is a homomorphism.
\end{remark}

\begin{defn}\label{compvoldefn} Let $M$ be an oriented tame $3$--manifold and let $\rho$ be a boundary-parabolic representation. The \emph{complex volume} of $\rho$ is defined by 
\begin{equation}\label{genvol} i(\Vol(\rho)+i\CS(\rho))=\widehat L\circ \Psi (\alpha),
\end{equation} where $\alpha$ is the fundamental class of some decoration of $\rho$.
\end{defn}

By Remark \ref{semi} and Theorem \ref{Lhat}, this definition agrees with the usual definition if $M$ is a hyperbolic manifold and $\rho$ is the geometric representation.

\begin{prop}\label{changes} The complex volume is unchanged when $\rho$ is changed by conjugation by an element in $\PSL(2,\C)$. If we change the orientation of $M$, the complex volume changes sign. If we change $\rho$ by composing with the involution on $\PSL(2,\C)$ given by complex conjugation, $\Vol(\rho)$ changes sign while $\CS(\rho)$ is fixed. 
\end{prop}
\begin{proof} The first statement follows from the fact that the associated $(G,P)$--cocycle only depends on the conjugation class of $\rho$. The second statement is obvious; all that changes is the signs of the $\epsilon_i$'s in \eqref{fundclassrep}. A change of $\rho$ by complex conjugation corresponds to changing the $(G,P)$--cocycle by complex conjugation. This changes all the log-parameters by complex conjugation, and since 
\[\widehat L([\bar z,-p,-q])=\overline{\widehat L([z,p,q])},\]
the third statement follows.
\end{proof}


\begin{ex}
We continue the study of the $5_2$ knot complement from Example~\ref{knotex}. The fundamental class of a representation conjugate to a representation in $P$ is clearly trivial, and since the $5_2$ knot is a twist knot, we know from Remark~\ref{twobridge} that all other boundary-parabolic representations are Galois conjugates of the geometric one. 

We first compute the complex volume of the geometric representation. Recall that the cross-ratios are given by $u=x^2$, $v=x^2-x+1$ and $w=v$, where $x$ is a root of $x^3-x^2+1$. For the geometric solution, the approximate values are
\[x=0.8774+0.7448i,\quad u=0.2150+1.3071i,\quad v=w=0.3376+0.5622i.\]
The flattenings are given by \eqref{logpar}. Using this together with \eqref{matrices} and the fact that $\sqrt{1-x}=1-v$ and $\sqrt{x^2-x-1}=u$ for the value of $x$ corresponding to the geometric solution, we see that the flattening $(w_0,w_1,w_2)$ for the first simplex is given by 
\begin{gather*}\label{firstlogpars}
\begin{aligned} w_0=&\Log(1-v)+\Log(u)-\Log(1)-\Log(1-v)=\Log(u)\\
w_1=&\Log(1)+\Log(1-v)-\Log(1)-\Log(u)=-\Log(1-u)-\pi i\\
w_2=&\Log(1)+\Log(u)-\Log(1-v)-\Log(u)=-w_0-w_1.
\end{aligned}\end{gather*} 
Omitting the calculation of $w_2=-w_0-w_1$, the flattenings of the second and third simplex are given as follows:

\begin{gather*}\label{secondlogpars}
\begin{aligned}
w_0=&\Log(1-v)+\Log(1)-\Log(u)-\Log(1)=\Log(v)-\pi i\\
w_1=&\Log(u)+\Log(1)-\Log(1-v)-\Log(u)=-\Log(1-v)
\end{aligned}
\end{gather*}

\begin{gather*}\label{thirdlogpars}
\begin{aligned}
w_0=&\Log(1-v)+\Log(1-v)-\Log(1-v)-\Log(u)=\Log(w)-\pi i\\
w_1=&\Log(1-v)+\Log(u)-\Log(1)-\Log(1)=-\Log(1-w)
\end{aligned}
\end{gather*}
Hence, we obtain that the image of the fundamental class in $\widehat \B(\C)$ is given by
\[[u;0,-1]+[v;-1,0]+[w;-1,0]\in\widehat \B(\C).\]
From this we can calculate the complex volume of the $5_2$ knot complement to be
\begin{equation}\label{georepknot}2.828122088330783\ldots+i\,3.024128376509301\ldots\in \C/i\pi^2\Z.\end{equation}

It is well known that the trace field (see Example~\ref{galconj}) of a link complement is equal to the field generated by the cross-ratio parameters. By Example~\ref{knotex}, the trace field of the $5_2$ knot complement thus equals $\Q(x)$. Since $x$ satisfies $x^3-x^2+1=0$, this field has degree $3$, and therefore has one real embedding and two complex conjugate embeddings. By Proposition \ref{changes}, the complex volume of the complex conjugate of the geometric representation is given by \eqref{georepknot}, but with the real part having the opposite sign. Let us compute the complex volume of the real Galois conjugate. The cross-ratio parameters are still given by $u=x^2$, $v=w=x^2-x+1$, but now $x$ is the real solution to $x^3-x^2+1=0$.
By the exact same method as above, but with $1-v$ replaced by $v-1$, since $\sqrt{1-x}=v-1$ for the real solution, we obtain that the element in $\widehat \B(\C)$ is
\[[u;0,0]+[v;0,1]+[w;0,1]\in \widehat\B(\C),\]
which gives a complex volume of
\[-i\,1.1134545524739240\ldots \in \C/i\pi^2\Z.\]

\end{ex}

\begin{remark} Note that essentially the same formula applies to calculate the log-parameters for both the geometric representation and its Galois conjugates. If we had used Neumann's formula we would have had to solve a new set of linear equations for each Galois conjugate. 
\end{remark}

\begin{remark} The $5_2$ knot complement is listed as $m015$ in the Snap census and Snap computes its Chern--Simons invariant to be $-3.024...$ (mod $\pi^2)$, which has the opposite sign of our result \eqref{georepknot}. This is because the census manifold is the mirror image of the standard $5_2$ knot complement from the Rolfsen table.
\end{remark}

\section{Lifts of boundary-parabolic representations}\label{SLreps}
This section is devoted to a discussion of representations in $\SL(2,\C)$. 
We shall see that a boundary-parabolic $\SL(2,\C)$--representation
has a fundamental class in $H_3(\SL(2,\C))$, and that a hyperbolic manifold with a spin structure has a fundamental class  in $H_3(\SL(2,\C))$ defined modulo $2$--torsion. The $2$--torsion ambiguity has the interesting consequence that a large class of cusped hyperbolic manifolds don't have ideal triangulations admitting even flattenings.

By the methods of Section \ref{fundclasssec}, a decorated boundary-parabolic 
$\SL(2,\C)$--repre\-sentation determines a fundamental class in
$H_3(\SL(2,\C),P)$. We wish to lift the map $\Psi\colon
H_3(\PSL(2,\C),P)\to \widehat \B(\C)$ to a map defined on
$H_3(\SL(2,\C),P)$ and taking values in the \emph{more} extended Bloch
group. The more extended Bloch group is defined as in Definition \ref{Bhat} but without including the transfer relation, and requiring that the integers $p$ and $q$
in Definition \ref{logpardefn} be even. We will refer to such flattenings as \emph{even} flattenings. The more extended Bloch group is shown in Goette--Zickert~\cite{GZ} to be isomorphic to $H_3(\SL(2,\C))$. We will therefore in this section denote the two versions of the extended Bloch group by $\widehat\B_{\SL}(\C)$ and $\widehat \B_{\PSL}(\C)$, respectively.

In Dupont--Zickert~\cite{DupontZickert} we studied a complex $C_*^{h\neq}(\C^2)$ generated in dimension $n$ by $n$--tuples of complex vectors in $\C^2$ in general position.
Using the simple fact that $\SL(2,\C)/P$ equals $\C^2-\{0\}$, it is not difficult to see that the cokernel of the map $C_2^{h\neq}(\C^2)\to C_1^{h\neq}(\C^2)$ is equal to the kernel of the augmentation map $C_0(\SL(2,\C)/P)\to \Z$. By Theorem~\ref{relhom} we have a canonical isomorphism 
\[H_3(\SL(2,\C),P)\cong H_3(C_*^{h\neq}(\C^2)\otimes_{\Z[\SL(2,\C)]}\Z).\]
An explicit formula for this isomorphism is given by
\begin{gather}\bar C_n(\SL(2,\C),P)\to C_n^{h\neq}(\C^2)\nonumber\\
\{g^{ij}\}\mapsto (v_0,\dots,v_n),\end{gather} where $v_i=g^{ij}\infty$, which is independent of $j$.

We can now define the desired lift of $\Psi$ as the composition of the isomorphism above with the explicit map
\[H_3(C_*^{h\neq}(\C^2)\otimes_{\Z[\SL(2,\C)]}\Z)\to \widehat \B_{\SL}(\C)\] constructed in Dupont--Zickert~\cite{DupontZickert}. The map $\widehat L$ defined in Theorem~\ref{Lhat} is shown in Goette--Zickert~\cite{GZ} to have a lift to $\widehat \B_{\SL}(\C)$ taking values in $\C/4\pi^2\Z$.
Summarizing the above, we have a commutative diagram with exact columns.
\[\xymatrix{{\Z/4\Z}\ar@{_{(}->}[d]\ar@2{-}[r]&{\Z/4\Z}\ar@{_{(}->}[d]\ar@2{-}[r]&{\Z/4\Z}\ar@{_{(}->}[d]\\
{H_3(\SL(2,\C),P)}\ar[r]^-\Psi\ar@{>>}[d]&{\widehat \B_{\SL}(\C)}\ar@{>>}[d]\ar[r]^{\widehat L}&{\C/4\pi^2\Z}\ar@{>>}[d]\\
{H_3(\PSL(2,\C),P)}\ar[r]^-\Psi&{\widehat \B_{\PSL}(\C)}\ar[r]^{\widehat L}&{\C/\pi^2\Z}}\]
The commutativity of the lower left square follows from the fact that for even flattenings, a semi-strong flattening is always a strong flattening. The proof of this follows the proof of Neumann~\cite[Proposition 5.3]{Neumann} word by word.

We can now use the above diagram to define the complex volume of a boundary-parabolic $\SL(2,\C)$--representation as an element in $\C/4i\pi^2\Z$.

\subsection{Spin structures and even flattenings}

A spin structure of a hyperbolic manifold $M$ is equivalent to a lift of the geometric representation to
$\SL(2,\C)$. If $M$ is closed it thus defines a fundamental class in $H_3(\SL(2,\C))\cong\widehat \B_{\SL}(\C)$, and a complex volume in $\C/4i\pi^2\Z$. If $M$ has cusps, lifts of the geometric representation are not boundary-parabolic. They are only $(\SL(2,\C),\pm P)$--representations. In fact, the proposition below, see e.g. Calegari~\cite[Corollary 2.4]{Calegari}, gives a concrete obstruction to defining a fundamental class in $H_3(\SL(2,\C))$ of a cusped hyperbolic manifold with a spin structure.
\begin{prop}\label{obstruct} Let $M$ be a cusped hyperbolic manifold. Any lift of the geometric representation to $\SL(2,\C)$ maps any curve bounding a $2$--sided incompressible surface to an element of $\SL(2,\C)$ with trace $-2$. 
\end{prop}

\begin{remark}It is not difficult to check that the map \begin{equation*}\label{spinstructure}
H_3(\SL(2,\C),P)\to H_3(\SL(2,\C),\pm P)
\end{equation*} is surjective with kernel of order $2$. This implies that a hyperbolic manifold with a spin structure does have a fundamental class in $H_3(\SL(2,\C))$ modulo $2$--torsion.
\end{remark}

\begin{remark}\label{neverSL}
Recall that a decorated ideal triangulation of a cusped hyperbolic manifold $M$ naturally gives rise to a $(\PSL(2,\C),P)$--cocycle on $\bar M$ representing the geometric representation.
As explained in Remark \ref{SL}, the labelings can be regarded as elements in $\SL(2,\C)$, but by Proposition \ref{obstruct} above, the $(\PSL(2,\C),P)$--cocycle is never an $(\SL(2,\C),P)$--cocycle.  
\end{remark}
Proposition \ref{obstruct} has the following interesting consequence.
\begin{thm}
Let $M$ be a cusped hyperbolic $3$--manifold satisfying

\begin{equation}\label{noeven}
\Ker\big(H^1(M;\Z/2\Z)\to H^1(\partial \bar M;\Z/2\Z)\big)=0.
\end{equation}

There is no ideal triangulation of $M$ admitting an even, strong flattening.
\end{thm}
\begin{proof}
As in \eqref{flatj}, we can regard any flattening as an element in $J\otimes \C$. Pick a decorated ideal triangulation of $M$, and let $w\in J\otimes \C$ be the flattening given by the decoration. Since $M$ satisfies \eqref{noeven}, it follows from Remark \ref{parityhomo} that $w$ is a strong flattening. By Lemma \ref{betaeqw}, $w=-\epsilon_i \beta (l)$, where $l\in C_1\otimes \C$ is given by the $\Log(c)$--parameters. It is a simple consequence of Remark~\ref{neverSL} and Remark~\ref{pqeven} that $w$ is \emph{not} an even flattening. By Neumann \cite[Lemma 9.3]{Neumann}, any two strong flattenings differ by an element in $C_1\otimes \pi i\Z$. Hence, any strong flattening $w'$ is given as $-\epsilon_i\beta(l+x)$, where $x$ is an element in $C_1\otimes \pi i \Z$. Let $E$ be a one-cell regarded as an element of $C_1=C_1\otimes \Z$. Recall that $E$ has a labeling \[g(E)=\begin{pmatrix}0&-c^{-1}\\c&0\end{pmatrix}\in\SL(2,\C),\] with $\Log(c)$--parameter $\Log(c)$. Note that adding $\pi i E\in C_1\otimes \pi i\Z$ to $l$ has the same effect on the parity of the corresponding flattening as changing the sign of $g(E)$. Now if $w'$ were an even flattening, this would imply that it would be possible to obtain an even flattening by changing signs of some of the labelings of long edges. By Remark \ref{pqeven}, this would imply that the new labelings (after a global sign change if necessary) would constitute an $(\SL(2,\C),P)$--cocycle. This is impossible by Proposition~\ref{obstruct}.
\end{proof}

\begin{remark} It seems worth mentioning that the property \eqref{noeven} also implies that the trace field of $M$ is equal to the invariant trace field of $M$. 

\end{remark}


\bibliographystyle{plain}
\bibliography{mybib}\end{document}